\documentclass{amsart}
\usepackage[T1]{fontenc}
\usepackage{amssymb}
\usepackage[hyperref,nonotes]{degt}

\def\2{\color{red}}
\let\2\relax

\def\itemrefform#1{\rom{(#1)}}

\makeatletter
\def\HyPsd@CatcodeWarning#1{}
\makeatother

\addtheorem{Claim}{claim}

\def\dual{^\vee}
\def\even{\sb\text{\rm even}}
\def\smod{\pdfstr\textasciitilde{^\sim}}
\def\smod{\pdfstr{\sp\#}{^\sharp}}
\def\units{^\times}
\def\sdif{\mathbin\vartriangle}
\let\ul\underline
\let\Leech\Lambda

\def\bL{\bold{L}}
\def\bH{\bold{H}}
\def\vv#1{\bold1_{#1}}
\def\be#1{\bar e_{#1}}
\def\be#1{[#1]}
\def\bo{\bar o}
\def\IS{\Cal I}
\def\GD{\Delta}
\def\Fn{\operatorname{Fn}}
\def\mFn{\underline\Fn}
\def\spn{\operatorname{span}}
\def\spnZ{\spn_\Z}
\def\spnO{\spn_0}

\def\sat{\operatorname{sat}}
\def\ort{\operatorname{ort}_h}
\def\FF{\Bbb{F}_2}
\def\Bnd{\Cal{B}}

\def\CA{\Cal A}
\def\CC{\Cal C}
\def\CD{\Cal D}
\def\CP{\Cal P}
\def\fS{\frak S}

\def\fC{\frak C}
\def\sS{\Cal S}

\def\sQ{\Cal Q}

\def\B(#1,#2){N_{#1}(#2)}
\def\Bs(#1,#2){N^*_{#1}(#2)}
\def\bN{N}
\def\N{N}
\def\RG{\GROUP{R}}
\def\orb{\frak o}
\def\Orb{\frak O}
\def\cluster{\frak c}
\def\borb{\bar\orb}
\def\cnt{c}
\def\bnd{b}
\def\bb{\bar\bnd}
\def\rt{\operatorname{\frak{rt}}}
\def\Sym{\operatorname{Sym}}

\let\set\Omega
\def\indA(#1,#2){^{#1}_{#2}}
\def\indA(#1,#2){_{#1,#2}}
\def\bbA{\CA\indA}

\def\bbB(#1,#2){\binom{#2}{#1}}
\def\bbC{\CC\indA}
\let\bbB\bbC
\def\indD(#1){_{#1}}
\def\bbD{\CD\indD}

\def\bbP{\CP\indD}
\def\bbT{\Cal T\indD}

\def\Golay{\Cal G}

\let\cl\ul
\def\orbit#1{\ifcase#1\or\or2_\circ\or\or4_\bullet\fi}

\def\hk{h|_k}
\def\orbk{\orb|_k}
\def\prim#1{#1^*}

\newcounter{subeq}[equation]

\def\makesubeq{\refstepcounter{subeq}\hbox{(\alph{subeq})\enspace}}
\def\subeq{\relax\let\label\savedlabel
 \ifmmode\makesubeq\else\par\strut\makesubeq\fi\ignorespaces}

\let\plussign+
\let\minussign-
\let\eqsign==
\let\black\bullet
\let\white\circ
\let\void\cdot
\def\ministrut{\vrule width0pt height8pt depth2pt}
\def\minibox{\hbox to 1em}
\def\minisize{\scriptstyle}
\def\LET{\let\\\cr}
\newcount\minicount
\chardef\minimax=100
{\catcode`\*\active\catcode`\.\active\catcode`\!\active
\catcode`\+\active\catcode`\-\active\catcode`\=\active
\gdef\miniset{\vcenter\bgroup
 \catcode`\*\active\def*{&\black}%
 \catcode`\.\active\def.{&\void}%
 \catcode`\!\active\def!{&\white}%
 \catcode`\+\active\def+{&\plussign}%
 \catcode`\-\active\def-{&\minussign}%
 \catcode`\=\active\def={&\eqsign}%
 \LET\offinterlineskip
 \ialign\bgroup\ministrut\hss$\minisize##$\global\minicount0\relax&&
  \global\advance\minicount by1\relax
  \ifnum\minicount<\minimax\minibox{\hss$\minisize##$\hss}\fi\cr}}
\def\endminiset{\crcr\egroup\egroup}

\def\Tsize{\scriptstyle}
\def\Tsize{}
\def\minilist{\bgroup
 \setbox0\hbox{$\Tsize[00,0,00]$}\edef\Isize{\the\wd0}%
 \setbox0\hbox{$\Tsize(0,0)$}\edef\IIsize{\the\wd0}%
 \def\I##1{\hbox to\Isize{\hss$\Tsize##1$}}%
 \def\II##1{\hbox to\IIsize{\hss$\Tsize##1$}}%
 \def\T##1{\quad\vcenter{\offinterlineskip\halign{\strut
  \I{####}&\ \II{####}\cr##1\crcr}}\kern-20pt}%
 \def\setsize##1{\kern-20pt \ls|\fC|=##1\:\quad}%
 \def\sym##1{}%
 \def\ministrut{\vrule width0pt height7pt depth1pt}%
 \def\minibox{\hbox to 8pt}%
}
\def\endminilist{\egroup}

\let\savedlabel\label
\allowdisplaybreaks

\relax
\def\tabrefs#1{\cite[#1]{degt:800.tab}}
\def\tabref#1{\tabrefs{\autoref*{#1}}}
\title{Conics on smooth quartic surfaces}

\author{Alex Degtyarev}

\address{%
Department of Mathematics\\
Bilkent University\\
06800 Ankara, TURKEY}

\email{
degt@fen.bilkent.edu.tr}

\thanks{%
The author was partially supported by the T\"{U}B\DOTaccent{I}TAK grant 123F111%
}



\dedicatory{to Ay\c{s}e}

\begin{document}

\begin{abstract}
We prove that the maximal number of conics,
\latin{a priori} irreducible of reducible,
on a smooth spatial quartic
surface is~$800$, realized by a unique quartic. We also classify quartics
with many (at least $720$) conics. The maximal number of real conics on a
real quartics is between $656$ and $718$.
\end{abstract}

%% Temporary
%\include{../800/800}

\maketitle

\section{Introduction}\label{S.intro}

Unless stated otherwise, all algebraic varieties in this paper are
over~$\C$.

\subsection{Principal results}\label{s.results}
Following~\cite{Bauer:conics}, denote by $\B(2n,d)$ the maximal number of
smooth rational degree~$d$ curves that a
smooth $2n$-polarized $K3$-surface $X\subset\Cp{n+1}$ may contain.
The numbers $\B(2n,1)$ of lines, most notably $\B(4,1)$, have a long
history going back at least to F.~Schur~\cite{Schur:quartics} and currently are well
known, see
\cite{degt:lines,degt:sextics,DIS,rams.schuett,Schur:quartics,Segre,Veniani}
and further references therein.

{\2However,\mnote{\2new par, extended}
even the very next case, \viz. the numbers $\B(2n,2)$ of conics, is still wide
open.
It is fairly obvious~\cite{degt:lines,rams.schutt:24curves}
that, for $n\gg d$, the degree up to~$d$ smooth rational curves on a
$2n$-polarized $K3$-surface~$X$ are either linearly independent in $\NS(X)$
(so that their number does not exceed~$19$) or among the fiber components of an
elliptic pencil on~$X$; hence, $\B(2n,d)\le24$ for $n\gg d$.
In particular, $\B(2n,2)=24$ for $n\gg0$, see~\cite{Rams.Schutt:12A1};
the other values
of~$d$ are discussed in~\cite{degt.Rams:periodicity}.
(Alternatively, this bound follows from Miyaoka~\cite{Miyaoka}, where it
extends to singular irreducible rational curves as well.)
Thus, it is the low polarization degrees $2n$ that are of interest, and here
we have but
a few sporadic examples and conjectures
\cite{Barth.Bauer:conics,Bauer:conics,degt:800,degt:conics,degt:8Kummer,degt:4Kummer}.}

The principal result of the present paper is \autoref{th.main} below. In the spirit
of \cite{degt:lines,DIS} and other papers based on the global Torelli theorem
and lattice theory, we do not just prove a
%sharp upper
bound; {\2instead, we} classify
all large configurations of conics. Note also that, in spite of the
definition of $\B(2n,2)$,
\latin{a priori} we do not distinguish between irreducible and reducible
conics; thus, a \emph{conic} is either a smooth (planar) curve of
projective degree~$2$ or a pair of distinct intersecting lines.
Throughout the paper, if reducible conics are present, we describe the total
number of conics as
\[*
(\text{irreducible})+(\text{reducible})=(\text{total}).
\]
{\2
This\mnote{\2$\Bs(2n,2)$ moved here and extended}
approach is motivated by the
conjecture~\cite{degt:conics,degt:8Kummer,degt:4Kummer} stating that the
presence of
reducible conics does not affect the ultimate bound $\B(2n,2)$;
furthermore, there exists a lower bound $\Bs(2n,2)<\B(2n,2)$ such that
any smooth $K3$-surface $X\subset\Cp{n+1}$ of degree $2n$
with more than $\Bs(2n,2)$ conics (irreducible or reducible) has no lines,
so that all conics are irreducible. A similar bound
$\Bs(2n,1)<\B(2n,1)$ on the number of lines in the presence of exceptional
divisors is conjectured and substantiated
in \cite{degt.Rams:octics,degt.Rams:quartics}.
All these observations agree with
Miyaoka~\cite{Miyaoka}, where, in fact, a certain weighted combination of \emph{all}
curves of degree up to~$d$ is bounded, but the ``correct'' meaning of the
general $\Bs(2n,d)<\B(2n,d)$ is yet to be understood.
(Unlike $\Bs(2n,2)$, in $\Bs(2n,1)$ we count \emph{irreducible}
lines only; otherwise, the numbers may blow up.)}

\theorem[see \autoref{proof.main}]\label{th.main}
A smooth quartic
$X\subset\Cp3$ with at least $720$ conics
\rom({\2either} irreducible or reducible\rom) is projectively
equivalent to one of the following surfaces\rom:
\roster
\item\label{i.M20}
the \emph{$M_{20}$-quartic}~\eqref{eq.24A1.800.1}, see \cite{Mukai},
with $800$ irreducible conics,
\[*
z_0^4+z_1^4+z_2^4+z_3^4+12z_0z_1z_2z_3=0,\quad\text{or}
\]
\item\label{i.4A4}
a \emph{$4^2\frak{A}_4$-quartic}~\eqref{eq.24A1.736},
with $736$ irreducible conics \rom(\cf. \autoref{ad.family}\rom), or
\item\label{i.L2(7)}
one of the two \emph{$L_2(7)$-quartics}~\eqref{eq.24A1.728}
with $728$ irreducible conics, or
\item\label{i.T192}
the \emph{$T_{192}$-}, \latin{aka} \emph{Schur's quartic},
with $64$ lines and $144+576=720$ conics,
\[*
z_0(z_0^3-z_1^3)=z_2(z_2^3-z_3^3).
\]
\endroster
The quartics in items~\iref{i.M20}, \iref{i.4A4}, \iref{i.L2(7)} have no
lines\rom; hence, all conics are irreducible.
\endtheorem

Thus, \autoref{th.main} states that $\B(4,2)=800$,
as conjectured
in~\cite{degt:800,degt:4Kummer},\mnote{\2shortened; mostly moved up}
{\2and that the hypothetical bound $\Bs(4,2)=720$ is indeed well defined.
For the moment, the only other value known to be well defined is
$249\le\Bs(6,2)<261$ \vs. $\B(6,2)=285$, see~\cite{degt:conics}.
In a forthcoming paper, we will show that
$\Bs(8,2)=136$ \vs. $\B(8,2)=176$.}

The quartic in item~\iref{i.M20} was discovered (in conjunction with the
conic counting problem) in~\cite{degt:800}, upon which X.~Roulleau (private
communication) observed that the surface admits a faithful symplectic action
of~$M_{20}$ (hence the equation above) and as such it
was studied by C.~Bonnaf\'{e} and A.~Sarti~\cite{Bonnafe.Sarti}; later,
B.~Naskr\k{e}cki~\cite{Naskrecki} found explicit equations of all $800$ conics.
Quartic~\iref{i.4A4} has previously appeared in \cite{degt:4Kummer},
whereas \iref{i.L2(7)} and~\iref{i.T192} appeared in \cite{degt:8Kummer},
as examples of quartics with many conics.
Of course, Schur's quartic~\iref{i.T192} has been known ever since
Schur~\cite{Schur:quartics}, and it is
%a truly remarkable surface
famous for quite a few other extremal properties.
Thus, it is the only smooth quartic
\roster*
\item
maximizing the number $\Bs(4,2)=720$ of conics in the presence of lines,
\item
maximizing the number of \emph{reducible} conics on a smooth quartic,
see~\cite{degt:4Kummer},
\item
admitting a faithful symplectic action of Mukai's group $T_{192}$,
see~\cite{degt:8Kummer},
\item
maximizing the number $\B(4,1)=64$ of lines, see~\cite{DIS},
\item
minimizing the determinant $\det T\ge48$ of the transcendental lattice of a
singular $K3$-surface admitting a smooth quartic model,
see~\cite{degt:singular.K3}.
\endroster

\remark
The notation/terminology according to the group
$\Sym_hX\subset\Aut_hX$ of symplectic projective automorphisms
of~$X$ used in \autoref{th.main} is mainly due to
%the lack of imagination.
{\2lack of a better choice.}
This is justified for $M_{20}$ and $T_{192}$, as quartics admitting symplectic
actions of these groups are indeed unique, see
\cite{Bonnafe.Sarti,degt:4Kummer} and~\cite{degt:8Kummer}, respectively.
However,
$4^2\frak{A}_4$-quartics constitute a $1$-parameter family
(see \autoref{ad.family}), and only one of
them has $736$ conics. Similarly, there are two $\PGL(4,\C)$-classes of
$L_2(7)$-quartics (see,
essentially,~\cite{Hashimoto}): one
%of them
is item~\iref{i.L2(7)} in
\autoref{th.main} whereas the other has neither lines nor conics
as the polarization~$h$
%splits as a
{\2generates an orthogonal}
direct summand of the
N\'{e}ron--Severi lattice~$\NS(X)$.
\endremark

\subsection{Further observations and examples}\label{s.further}
Our principal results are proved by a computation and, in order to produce
some interesting examples, we saved and analyzed all large (at least $600$
conics) configurations encountered in the course of the proof. The
total conic counts observed are
\[*
\hyperref[i.M20]{800},
\hyperref[i.4A4]{736},
\hyperref[i.L2(7)]{728},
\hyperref[i.T192]{720}, 704,
\hyperref[eq.24A1.680]{680}, 664,
\hyperref[eq.24A1.660]{660},
\hyperref[ad.real]{656},
\hyperref[ad.family]{640}, 636, 628, 624, 622, 620, 616,
\hyperref[ad.family]{608}, 600,
\]
some of which are represented by multiple configurations/quartics.
We do not assert that this list is complete, but it must be close to such; in
particular, it does contain all Barth--Bauer quartics with at least $600$ conics,
see~\cite{degt:4Kummer}.

In particular, one of the examples found beats the record for the number of
\emph{real} conics on a \emph{real} quartic set
in \cite{degt:4Kummer}, see \autoref{ad.real}; the next
%known
{\2largest}
configuration
{\2(among those known)}
is $188 + 448 = 636$ real conics on $Y_{56}$ in~\cite{DIS}, see
\autoref{rem.reducible}.

\addendum[see \autoref{proof.real}]\label{ad.real}
There exists a real smooth quartic with $656$ real conics,
see \eqref{eq.24A1.656}\rom;
moreover, all conics are irreducible and have non-empty real part.
Hence, we have the bounds
\[*
656\le\B(4,2;\R)\le718
\]
on the maximal number of real conics on a real smooth quartic.
\endaddendum

Yet another example is the largest known configuration in a $1$-parameter
family. This family connects quartics~\iref{i.M20}
and~\iref{i.4A4} in \autoref{th.main}.
{\2Here and henceforth,\mnote{\2explained}
the terms \emph{equiconical} (\emph{equilinear}, \etc.)\ are supposed to be
self-explanatory: it suffices to assume that the \emph{numbers} of
irreducible/reducible conics remain constant within the family;
then, so do respective dual adjacency graphs, \cf.~\eqref{eq.conics}.}

\addendum[see \autoref{proof.family}]\label{ad.family}
There exists an equiconical $1$-parameter family~$\Cal X$ of \rom(generically\rom)
$4^2\frak{A}_4$-quartics with $608$ irreducible conics, see
\eqref{eq.24A1.608}.
The closure $\bar{\Cal X}$ of this
family contains
quartics~\iref{i.M20}
and~\iref{i.4A4} in \autoref{th.main}, as well as the Barth--Bauer quartic
with $640$ conics, see~\cite{degt:4Kummer} and~\eqref{eq.24A1.640}.
\endaddendum

The previously known record for both the number of real conics and the
number of conics in a family was $560$ (see
\cite{degt:8Kummer}; same quartic);
{\2for some reason, the  quartic $Y_{56}$ known
since~\cite{DIS} has never been analyzed.
We close this gap in \autoref{rem.reducible}.}

A few other examples are discussed in \autoref{s.remarks}.

\subsection{Idea of the proof}\label{s.idea}
Similar to \cite{degt:conics,degt:sextics,DIO}, we establish that, regarded as an
abstract graph, the configuration of conics (irreducible or reducible) on a
smooth quartic surface can be realized as a certain special set of square~$4$
vectors in a $4$-polarized Niemeier
lattice, \ie, positive definite even unimodular lattice of rank~$24$.
There are but $24$ Niemeier lattices, all of which are well known, see
\cite{Niemeier,Conway.Sloane}; this already implies that there is indeed a
uniform bound on the number of conics, a fact that is not immediately obvious
\latin{a priori}.

The principal novelty is the fact that, due to the much larger numbers
involved, in order to make the computation feasible we have to shift the
paradigm and deal with sublattices (or even rational subspaces) rather
%that just subsets.
{\2than just admissible sets of lines in the sense of \cite{degt:sextics}.}
This shift lets us revise and substantially refine the
combinatorial estimates on the number of conics, immediately resulting in the
bound of $1736$, which is much better than the previously known
$5016$ (see [2], with a reference to S.~A.~Str{\o}mme). Reducing this further
down to $718$ (with a few exceptions stated in \autoref{th.main}), still
computer aided, requires a drastic revision of all parts of the algorithm and
the underlying mathematics: in addition to working with subspaces rather than
subsets, we employ the full orthogonal group rather than just reflections,
treat more carefully iterated index~$2$ subgroups in \autoref{s.subgroups},
\etc.

\subsection{Contents of the paper}\label{s.contents}
In \autoref{S.Niemeier} we reduce the conic counting problem to the study of
the so-called \emph{geometric} sets of square~$4$ vectors in $4$-polarized Niemeier
lattices.
The precise conditions are formalized in \autoref{s.Niemeier}, and
in the rest of \autoref{S.combinatorial} we deal with the combinatorial
background necessary for \autoref{S.many.roots}, where we rule out
$18$ out of the $24$ Niemeier lattices. Unlike a few earlier
papers~\cite{degt:sextics,DIO}, here we consider primitive sublattices, hence
so-called \emph{saturated} sets only; this fact lets us refine the concept of
admissibility and improve the combinatorial estimates. Our principal goal is
making this part of the proof as human comprehendible as possible;
several examples are worked out in full detail.

In \autoref{S.few}, we treat the remaining five lattices rationally generated
by roots, and this part of the proof, based essentially on \autoref{lem.chain},
is heavily computer aided. To save space and keep the code separate from the
underlying mathematical ideas, we have moved the details to the companion
text~\cite{degt:800.tab}.
Here, we merely state the result, outline the proof, and list exceptional or
otherwise interesting configurations of conics
in a form from which they can easily be reconstructed.

The Leech lattice is treated in a similar manner in \autoref{S.Leech}. The
heart of our
approach is \autoref{lem.subgroups}. We give a complete proof of this
statement since it is followed closely by the new version of the code, which
is much more efficient than that used in~\cite{degt:conics,DIO},
see \autoref{s.subgroups} and \tabref{S.Leech.tab}.
We discover that all extremal configurations of conics listed in
\autoref{th.main} can be embedded to the Leech lattice, as
described in \autoref{proof.Leech}.

In \autoref{S.proofs} we collect the output of the computation in
\autoref{S.many.roots}--\autoref{S.Leech} and complete the formal proof of the
principal results of the paper stated in the introduction.

All technical details and a plethora of further examples of configurations
of conics are found in the companion text~\cite{degt:800.tab}: it is
available electronically, {\2including the \GAP~\cite{GAP4.13} code,}
%from the author's web page or
as an ancillary file for this paper in
%the \texttt{arXiv}.
\arXiv{2407.00493}.
%\mnote{\2links}
%{\2The \GAP~\cite{GAP4.13} code is also available from the \texttt{arXiv}.}

\subsection{Acknowledgements}\label{s.acknowledgements}
This work was partially supported by the Scientific and Technological Research Council of
Turkey (T\"{U}B\DOTaccent{I}TAK), grant number 123F111. I am
grateful to T\"{U}B\DOTaccent{I}TAK for
their support.
The paper was conceived during my sabbatical stay at the
{\em Max-Planck-Institut f\"{u}r Mathematik}, Bonn; I extend my gratitude to
this institution and its friendly staff for the warm atmosphere and excellent
working conditions. My special thanks go to the MPIM's computer team: it is
their involvement that let me complete the computation.

Most algorithms used in the paper were implemented using \GAP~\cite{GAP4.13}.

\section{Embedding to a Niemeier lattice}\label{S.Niemeier}
In this section we exploit the fact that a smooth spatial quartic is a $K3$-surface
and reduce the conic counting problem to the study of large collections
of certain square~$4$
vectors in a Niemeier lattice.

\subsection{Smooth quartics as $K3$-surfaces}\label{s.K3}
Recall that a smooth quartic $X\subset\Cp3$ is a $K3$-surface; in particular,
there is an isometry
\[*
H_2(X)\cong\bL:=-2\bE_8\oplus3\bU.
\]
(Here and below, we use
Poincar\'{e} duality to identify $H^2(X)$ and $H_2(X)$ and,
\via\ the intersection index form,
regard the latter as
a unimodular integral lattice.) The N\'{e}ron--Severi lattice
$\NS(X)\subset H_2(X)$ is a primitive hyperbolic sublattice of
\emph{Picard rank} $\Gr:=\rank\NS(X)\le20$; this lattice is naturally
$4$-polarized by the class $h\in\NS(X)$, $h^2=4$, of hyperplane section.

The following fact is essentially contained in Saint-Donat~\cite{Saint-Donat}
(see also \cite{degt.Rams:octics} for an accurate restatement in terms of
homology classes rather than linear systems); regarding the homology,
we also refer to the surjectivity of the period map~\cite{Kulikov:periods}.

\theorem[Saint-Donat~\cite{Saint-Donat}]\label{th.Saint-Donat}
Let $S\ni h$ be a $4$-polarized hyperbolic lattice
and $S\into\bL$ a primitive isometric embedding. Then there is
a \emph{smooth} quartic $X\subset\Cp3$ and
an isometry
\[*
\bigl(\bL\supset S\ni h\bigr)\cong\bigl(H_2(X)\supset\NS(X)\ni h\bigr)
\]
%assume that $S$ admits a primitive
%embedding to $\bL$. Then $S\ni h$ is isomorphic to
%the N\'{e}ron--Sevri lattice $\NS(X)\ni h$ of a
%for a certain
%\emph{smooth} quartic $X\subset\Cp3$
if and only if there is no class $e\in S$ such that
\roster
\item\label{i.exceptional}
$e^2=-2$ and $e\cdot h=0$ \rom(\emph{exceptional divisor}\rom), or
\item\label{i.isotropic}
$e^2=0$ and $e\cdot h=2$ \rom(\emph{$2$-isotropic vector}\rom).
\done
\endroster
\endtheorem

Next,
we
{\2consider an arbitrary polarized hyperbolic lattice $S\ni h$\mnote{\2$S$
reintroduced}
and}
recall Vinberg's algorithm~\cite{Vinberg:polyhedron} for computing
the fundamental polyhedron of the group generated by reflections in the
hyperbolic space
\[*
\bigl\{x\in S\otimes\R\bigm|x^2>0\bigr\}/\R\units
\]
containing the class of~$h$ or, for short, just the \emph{fundamental
polyhedron of $S\ni h$}.
For each $n\in\NN$, introduce inductively the sets
\[*
\aligned
\GD_n(S,h)&:=\bigl\{x\in S\bigm|\text{$x^2=-2$ and $x\cdot h=n$}\bigr\},\\
\GD_n^\circ(S,h)&:=\bigl\{x\in\GD_n\bigm|
 \text{$x\cdot y\ge0$ for all $y\in\GD_m^\circ$, $m<n$}\bigr\}.
\endaligned
\]
%We use the fact that $X$ is assumed smooth and, by \autoref{th.Saint-Donat},
{\2Now,\mnote{\2assumptions made when needed}
assume that
$\GD_0^\circ=\GD_0=\varnothing$
(recall that we are interested in
\emph{smooth} surfaces, \cf. \autoref{th.Saint-Donat});}
this lets us avoid the subtlety
 with the choice of a Weyl chamber for the root system in
$h^\perp\subset S$ (\cf., \eg, \cite{degt.Rams:octics}).
Then, the fundamental polyhedron of $S\ni h$ is the set
\[*
\GD(S,h):=\bigcup_{n=0}^\infty\GD_n^\circ(S,h).
\]
(Strictly speaking, the fundamental polyhedron is the intersection of the
half-spaces $x\cdot v\ge0$, $v\in\GD$, but, since we are interested in
the set of its
walls only, we abuse the language and refer to $\GD$ itself as the
fundamental polyhedron.)

Now, assume that $(S\ni h)=(\NS(X)\ni h)$ for a quartic $X\subset\Cp3$.
Then, according to \cite[\S\,8.1]{Huybrechts}, the map
$C\mapsto[C]\in S$ establishes a bijection between the set of smooth rational
curves on~$X$ and $\GD$. Taking into account the projective degree $[C]\cdot h$
of a curve and still assuming that $X$ is smooth,
we have canonical bijections
\[
\alignedat2
\Fn_1(X,h)&:=\{\text{lines on~$X$}\}&&=\GD_1^\circ(S,h)=\GD_1(S,h),\\
\Fn_2^\circ(X,h)&:=\{\text{irreducible conics on~$X$}\}&&=\GD_2^\circ(S,h),\\
\Fn_2(X,h)&:=\{\text{all conics on~$X$}\}&&=\GD_2(S,h).\\
\endalignedat
\label{eq.conics}
\]
All three sets, as well as the
\emph{graph $\Fn_*(X,h):={\Fn_1}\cup{\Fn_2^\circ}$ of lines and conics}, are
regarded as (multi-)graphs, with two vertices represented by curves $C_1$, $C_2$
connected by an edge of multiplicity $C_1\cdot C_2$
(no edge if $C_1\cdot C_2=0$). It is easily seen that all multiplicities are
non-negative, even if both $C_1,C_2\in\Fn_2(X,h)$ are (distinct) reducible conics.
The graph $\Fn_*$ is also \emph{colored} according to the projective degree
of its vertices.

The {\2three} Fano graphs above can be defined for any
$4$-polarized hyperbolic lattice $S\ni h$ provided that $S$ does not contain
a vector~$e$ as in \autoref{th.Saint-Donat}\iref{i.exceptional}
or~\iref{i.isotropic}: it is this condition that keeps the multiplicities
non-negative, \cf. \autoref{lem.intr} below.

\subsection{The modified N\'{e}ron--Severi lattice}\label{s.S*}
Start with a $4$-polarized hyperbolic lattice $S\ni h$ and
%define its
%\emph{$h$-even sublattice}
consider the sublattice
\[*
S\even :=\bigl\{x\in S\bigm|x\cdot h=0\bmod2\bigr\}.
\]
Clearly, $S\even=S$ or $[S:S\even]=2$.
This lattice is still hyperbolic and
{\2$4$-polarized (by the same vector $h\in S\even$).
Besides, it is\mnote{\2rearranged}}
%it is \emph{$h$-even} in the sense that
%it is $4$-polarized (by the same vector $h\in S\even$) and
\[\label{eq.h-even}
\text{\emph{$h$-even}}\:\quad
\text{$h\in2S\even\dual$, \ie, $x\cdot h=0\bmod2$ for each $x\in S\even$}.
\]
Now, given an $h$-even $4$-polarized lattice $S\ni h$, we define $(S\ni h)\smod$
to be the same pair $S\ni h$ with the bilinear form on~$S$ modified \via
\[*
x\otimes y\mapsto\tfrac12(x\cdot h)(y\cdot h)-(x\cdot y).
\]
Since this operation makes sense for $h$-even lattices only, we will use the
shorthand $(S\ni h)\smod:=(S\even\ni h)\smod$ in the general case (\cf. also
the alternative construction in \autoref{s.embedding} below, where the passage
to $S\even$ is automatic).

\lemma\label{lem.involutive}
The map $(S\ni h)\mapsto(S\ni h)\smod$ is an involutive operation on the set of
$h$-even $4$-polarized even lattices. A lattice $S$ is hyperbolic if and
only if $(S\ni h)\smod$ is positive definite.
\endlemma

\proof
For the last assertion, observe that the modification preserves
%the square
$h^2=4$ and reverses the form on $h^\perp$. All other statements
are straightforward.
\endproof

For a positive definite $4$-polarized $h$-even lattice $S\smod\ni h$
we modify the
notion of its \emph{Fano graph}, letting
\[*
\mFn(S\smod,h):=\bigl\{x\in S\smod\bigm|
 \text{\rm$x^2=4$ and $x\cdot h=2$}\bigr\}
\]
and connecting two vertices $x,y$ by an edge of multiplicity
$2-x\cdot y$. As above, this notion is mostly useful if $S\smod$ is root
free; then, all multiplicities are non-negative
due to the following lemma.

\lemma\label{lem.intr}
If a positive definite $4$-polarized lattice $S\ni h$ is root free, then, for
any two vectors $l_1,l_2\in\mFn(S,h)$ one has
\[*
l_1\cdot l_2=-2\ (\text{\rm iff $l_1+l_2=h$}),\ 0,\ 1,\ 2,
 \ \text{\rm or}\ 4\ (\text{\rm iff $l_1=l_2$}).
\]
\endlemma

\proof
Consider the sublattice $R:=\Z h+\Z l_1+\Z l_2$ and assert that $\det R\ge0$.
This results in $-2\le l_1\cdot l_2\le4$. In the two border cases, $\det R=0$
and the vectors are linearly dependent, as stated in the lemma. If
$l_1\cdot l_2=3$, then $l_1-l_2$ is a root, and if $l_1\cdot l_2=-1$, then
$h-l_1-l_2$ is a root.
\endproof

\proposition\label{prop.even}
Assume that a $4$-polarized hyperbolic lattice $S\ni h$ admits a primitive
embedding to $\bL$. Then\rom:
\roster
\item\label{i.root.free}
$S\ni h$ is isomorphic to the N\'{e}ron--Severi lattice $\NS(X)\ni h$ of a
smooth quartic $X\subset\Cp3$ if and only if the modified
lattice $(S,h)\smod$ is root free\rom;
\item\label{i.conics}
if $(S,h)\smod$ is root free, the graph $\Fn_2(X,h)$ of all conics of
any quartic~$X$ as in
item~\rom{\iref{i.root.free}} is canonically isomorphic to the Fano graph
$\mFn(S\smod,h)$.
%\[*
%\Fn\smod(S\even\smod,h):=\bigl\{x\in S\even\smod\bigm|
% \text{\rm$x^2=4$ and $x\cdot h=2$}\bigr\},
%\]
%where two vertices $x$, $y$ are connected by an edge of multiplicity
%$2-x\cdot y$.
\endroster
\endproposition

\proof
For the first statement, observe that any vector $e\in S$ as in
\autoref{th.Saint-Donat}\iref{i.exceptional} or~\iref{i.isotropic} would
survive to $S\even$ and give rise to a root $e\in S\smod$ such that
$e\cdot h=0$ or~$2$, respectively. Conversely, any such root in~$S\smod$ is
a prohibited vector in~$S$. On the other hand, since $S\smod$ is positive
definite, $\ls|e\cdot h|\le2$ for any root
$e\in S\smod$. Since $S\smod$ is also $h$-even, this leaves
$e\cdot h\in\{0,\pm2\}$, \ie, at least one of $\pm e$ is
necessarily a prohibited vector.

The second statement follows from~\eqref{eq.conics}
and the obvious observation that
\[*
\GD_2(S,h)=\GD_2(S\even,h)=\mFn(S\smod,h).
\qedhere
\]
\endproof

\subsection{Embedding $S\smod$ to a Niemeier lattice}\label{s.embedding}
The chain
\[*
(S\ni h)\mapsto(S\even\ni h)\mapsto(S\even\ni h)\smod
\]
transforming a hyperbolic $4$-polarized lattice to a positive definite
$h$-even one can alternatively be described as follows:
\roster
\item\label{N.step.1}
consider the orthogonal complement $h^\perp\subset S$; change the sign of
the form;
\item\label{N.step.2}
consider the orthogonal direct sum $S'=(-h^\perp)\oplus\Z h$, $h^2=4$;
\item\label{N.step.3}
define $S\smod$ as the extension of $S'$ \via\ all/any vector of the form
\[*
\tfrac12(2l-h)\oplus\tfrac12h\in S'\otimes\Q,
\]
where $l\in S$ is such that $l\cdot h=2\bmod4$.
\endroster
If $l\in S$ as in Step~\iref{N.step.3} does not exist, we merely leave
$S\smod:=S'$; this case is not interesting as one obviously has
$\GD_2(S,h)=\varnothing$. Otherwise, both $S\even\supset S'$
(as an abelian group) and
$S\smod\supset S'$ are extensions of index $2$, and
\[*
l\mapsto\tfrac12(2l-h)\oplus\tfrac12h,\quad l\in S,\quad l\cdot h=2\bmod4,
\]
is a canonical bijection between the cosets
$S\even\sminus S'$ and $S\smod\sminus S'$.

It follows that the positive definite $4$-polarized lattice $S\smod\ni h$
constructed in \autoref{s.S*} from a hyperbolic lattice $S\ni h$
is the lattice $S(S,h)$ considered in \cite{degt:conics}. Thus, we have the
following statement, based on Nikulin's criterion~\cite{Nikulin:forms}.

\proposition[see {\cite[Proposition 2.10]{degt:conics}}]\label{prop.Niemeier}
If a $4$-polarized hyperbolic lattice $S\ni h$ admits a primitive embedding
to $\bL$, then the modified lattice $S\smod$ admits a primitive embedding to at
least one of the $24$ Niemeier lattices.
\done
\endproposition

%\subsection{An alternative construction}\label{s.alt}
%{\2\mnote{to do?}
%Change $S'$ in \iref{N.step.2} to $(-h^\perp)\oplus\Z r\oplus\Z s$,
%$r^2=s^2=2$, and do the extension \via\
%$\tfrac12(2l-h)\oplus\tfrac12r\oplus\tfrac12s$.
%}

\section{Combinatorial bounds}\label{S.combinatorial}
In this section we discuss a few simple (and very rough) combinatorial bounds
that rule out (in \autoref{S.many.roots} below) the vast majority of Niemeier
lattices.
The few remaining ones will be treated in the subsequent sections.

\subsection{Niemeier lattices\noaux{ (see~\cite{Conway.Sloane,Niemeier})}}\label{s.Niemeier}
Recall that, with one exception (the Leech lattice, see \autoref{S.Leech}
below), a Niemeier lattice~$\bN$ is rationally generated by roots and
is determined up to isomorphism by its maximal root
%system
{\2sublattice}~$D$.
Therefore, we use the notation
\[
\bN:=\N(D)=\N\bigl({\textstyle\bigoplus_kD_k}\bigr),\quad k\in\Omega,
\label{eq.Niemeier}
\]
where $D_k$ are the indecomposable $\bA$--$\bD$--$\bE$ components of~$D$ and
$\Omega$ is the index set.
There are well-defined orthogonal projections $N\to D_k\dual$, $k\in\Omega$,
and we also fix the notation
$l=\sum_kl|_k$, $l|_k\in D_k\dual$, $k\in\Omega$, for the decomposition of a
vector $l\in\bN$. The vector~$l$ or its projection~$l|_k$ is called
\emph{integral} if $l\in D$ or $l|_k\in D_k$, respectively; otherwise, if
$l\in\bN\sminus D$ or $l|_k\in D_k\dual\sminus D_k$, they are called
\emph{rational}. If $l^2=4$, then each projection~$l|_k$ is either integral,
and then $l|_k^2\in\{0,2,4\}$, or a shortest vector in its discriminant class
$(l|_k\bmod D_k)\in\discr D_k$, as otherwise a shorter representative would
give rise to a root $l'\notin D$. Furthermore, if $l^2=4$ and at least one
projection $l|_k\ne0$ is integral, then so are~$l$ and all other projections.

We fix a polarization $h\in\bN$, $h^2=4$, and abbreviate
\[*
\Orb:=\Orb_h:=\mFn(\bN,h)=\bigr\{l\in\bN\bigm|\text{$l^2=4$, $l\cdot h=2$}\bigr\}.
\]
The elements of $\Orb_h$ are called \emph{conics}.
We have a fixed point free \emph{duality} involution
\[*
{*}\:\Orb_h\to\Orb_h,\quad l\mapsto l^*:=h-l.
\]
Geometrically, $l$ and~$l^*$ represent a complementary pair of conics
constituting a hyperplane section of the quartic.
A subset $\fC\subset\Orb_h$ is \emph{self-dual} if $\fC^*=\fC$.

\definition\label{def.admissible}
For
subsets $\fC\subset\sS\subset\mFn(\bN,h)$, we define the \emph{spans}
\[*
\spnZ\fC:=\Z\fC+\Z h\subset\bN,\quad
\spn\fC:=\bN\cap(\Q\fC+\Q h)\subset\bN
\]
and \emph{saturations}
\[*
\sat_\sS\fC:=\sS\cap\spn\fC\subset\sS,\quad
\sat\fC:=\Orb_h\cap\spn\fC\subset\Orb_h.
\]
A subset $\fC$ is called
\roster*
\item
\emph{saturated}
(resp.\ \emph{$\sS$-saturated}) if $\sat\fC=\fC$ (resp.\ $\sat_\sS\fC=\fC$);
%it is called
\item
\emph{admissible} if $\spn\fC$ is
$h$-even and root free.
\endroster
Note that a saturated set is automatically self-dual;
if $\sS$ is self-dual, then so are all $\sS$-saturated sets.
\enddefinition

\autoref{def.admissible} extends to any $4$-polarized positive
definite lattice $\bN\ni h$. By an implicit reference to
$\spn\fC$ we also extend to subsets $\fC\subset\Orb_h$ such lattice-theoretic
notions as rank, discriminant, \etc., so that, \eg, $\rank\fC:=\rank\spn\fC$.

\definition\label{def.geometric}
An admissible subset $\fC\subset\mFn(\bN,h)$ is called \emph{geometric} if
there is an isometry
\[*
\Gf\:(\spn\fC\ni h)\smod\into\bL
\]
such that the primitive hull of the image $\Im\Gf$ is either $\Im\Gf$ itself
or a certain $h$-odd index~$2$ overlattice $S\supset\Im\Gf$.
\enddefinition

Since, for a $K3$-surface~$X$, one has $\Gs_+H_2(X)=19$ and $\NS(X)$ is
hyperbolic,
%and the lattice $\NS(X)$ is hyperbolic,
an obvious necessary condition for a set $\fC\subset\Orb_h$ to be geometric
is
that
\[
\rank\fC\le20;
\label{eq.rank}
\]
this is an essential part of all algorithms, even if we choose to skip the
expensive \cite[Proposition 2.10]{degt:conics}, possibly followed
by the analysis of index~$2$ extensions.

Note that, since $\spn\fC$ in \autoref{def.geometric} is assumed root free,
both $(\spn\fC\ni h)\smod$ and any $h$-odd index~$2$ extension thereof are
automatically free of vectors~$e\in S$ as in
\autoref{th.Saint-Donat}\iref{i.exceptional} or~\iref{i.isotropic}.
Furthermore,
according to \autoref{lem.intr}, an admissible subset $\fC\subset\mFn(\bN,h)$
can be regarded as a graph. Thus, \autoref{prop.even} and
\autoref{prop.Niemeier} imply the following statement.

\proposition\label{prop.reduction}
An abstract graph~$\Gamma$ is realizable as the graph $\Fn_2(X,h)$ of conics of a
smooth quartic $X\subset\Cp3$ if and only if $\Gamma$ is graph-isomorphic to
a saturated geometric subset $\fC\subset\mFn(\bN,h)$ in a $4$-polarized Niemeier
lattice $\bN\ni h$.
\done
\endproposition

\subsection{Orbits and combinatorial orbits}\label{s.orbits}
We make extensive use of the following
groups:
\roster*
\item
$\OG(N)$, the full orthogonal group of~$N$,
\item
$\RG(N)\subset\OG(N)$, the subgroup generated by reflections
$r_e\:x\mapsto x-(x\cdot e)e$ against the roots $e\in D$,
\item
$\OG_h(N)$, the stabilizer of~$h$ in $\OG(N)$, and
\item
$\RG_h(N):=\OG_h(N)\cap\RG(N)$.
\endroster
Since $\RG(N)$ acts simply transitively on the set of its fundamental
polyhedra, we conclude that
\[
\text{$\RG_h(\bN)$ is generated by reflections
 against the roots orthogonal to~$h$}.
\label{eq.Rh}
\]
Besides, for any overlattice or inner product $\Q$-vector space
$B\supset\bN$,
\[
\text{the action of $\RG_h(\bN)\subset\RG(\bN)$ extends to~$B\supset N$ \via\ reflections}.
\label{eq.action}
\]
%Note also that the groups $\RG_h(\bN)\subset\RG(\bN)$ act, \via\ reflections,
%on any overlattice or inner product $\Q$-vector space $B\supset\bN$.

The orbits of the action on~$\Orb_h$ of the groups $\OG_h(N)$ and $\RG_h(N)$
are called \emph{orbits} and \emph{combinatorial orbits}, respectively;
typically, each orbit~$\borb$ splits into a number of combinatorial orbits
$\orb_1,\orb_2,\ldots$, all of the same cardinality. The duality~$*$ induces an
involution (no longer free) on the sets of orbits and combinatorial orbits.

We define the \emph{count}
$\cnt(\orb)$ and \emph{bound} $\bnd(\orb)$ of a combinatorial orbit~$\orb$
\via
\[*
\cnt(\orb):=\ls|\orb|,\qquad \bnd(\orb):=\max\ls|\fC|,
\]
where $\fC$ runs over all \emph{admissible} subsets $\fC\subset\orb$. It is
obvious that $\cnt(\orb)$ and $\bnd(\orb)$ are constant within each orbit and
that $\cnt(\orb^*)=\cnt(\orb)$, $\bnd(\orb^*)=\bnd(\orb)$.

These na\"\i ve bounds are extended to unions of combinatorial orbits, \eg,
orbits~$\borb$ or the whole set~$\Orb_h$, \emph{by additivity}:
\[*
\cnt(\orb_1\cup\ldots)=\cnt(\orb_1)+\ldots,\qquad
\bnd(\orb_1\cup\ldots):=\bnd(\orb_1)+\ldots.
\]

\subsection{The lattices $\bH\sb n$, $\bA\sb n$, and $\bD\sb n$}\label{s.HAD}
Given an index set $\IS:=\IS_n:=\{1,\ldots,n\}$, we
consider the (odd unimodular) Euclidean lattice
\[*
\bH_{n}:=\bigoplus_{i\in\IS}\Z e_i,
\quad e^2_i=1.
\]
When $\bH_n$ and, hence, $\IS$ are fixed, we denote by
$\bo:=\IS\sminus o$ the complement of $o\subset\IS$ and
let $\vv{o}:=\sum e_i\in\bH_n$, $i\in o$.
The notation
$r\sdif s:=(r\cup s)\sminus(r\cap s)$ stands for the symmetric difference of sets.
We have
\[
{\subeq\label{eq.H.1}
\vv{u}\cdot\vv{v}=\ls|u\cap v|,\qquad
\subeq\label{eq.H.2}
 \vv{u}^2-\vv{u}\cdot\vv{v}=\frac12\ls|u\sdif v|\quad\text{if $\ls|u|=\ls|v|$}.
}\label{eq.H.prod}
\]
It is this and a few similar relations below in this section that
explain the relevance of bounds~\eqref{eq.bA}, \eqref{eq.bD} in
\autoref{s.combinat} below.

The group $\OG(\bH_n)$ is generated by permutations of the basis vectors
(equivalently, reflections against the roots $e_i-e_j$) and
reflections $r_i\:a\mapsto a-2(a\cdot e_i)e_i$ against the generators~$e_i$.
All roots in~$\bH_n$ are of the form $\pm e_i\pm e_j$, $i,j\in\IS$, $i\ne j$.

The root lattice $\bA_n$ is defined as $\vv{\IS}^\perp\subset\bH_{n+1}$:
\[*
\textstyle
\bA_n=\bigl\{\sum_i\Ga_ie_i\in\bH_{n+1}\bigm|\sum_i\Ga_i=0\bigr\}.
\]
The discriminant group is $\Z/(n+1)$, and the shortest vectors in the
discriminant class $\cl{m}=\cl0,\ldots,\cl{n}$
(we use the notation of~\cite{Conway.Sloane}, underlined to avoid confusion)
are
\[*
\be{o}:=\frac{1}{n+1}\bigl(\ls|\bo|\vv{o}-\ls|o|\vv{\bo}\bigr),\quad
 \be{o}^2=\frac{\ls|o|\ls|\bo|}{n+1},\quad
 \text{for $o\subset\IS$, $\ls|o|=m$},
\]
so that $\be{\bo}=-\be{o}$ and
\[
{\subeq\label{eq.A.1}
\be{r}\cdot\be{s}=\ls|r\cap s|-\frac{\ls|r|\ls|s|}{n+1},\qquad
%\text{and, for $\ls|r|=\ls|s|$},\quad
\subeq\label{eq.A.2}
\be{r}^2-\be{r}\cdot\be{s}=\frac12\ls|r\sdif s|\quad\text{if $\ls|r|=\ls|s|$}.
}\label{eq.A.prod}
\]
The $\RG(\bA_n)$-orbits of integral vectors of square~$4$ or~$2$ are
\[*
\orbit4:=\{\vv{u}-\vv{v}\},\quad\orbit2:=\{\vv{u}-\vv{v}\}\quad
\text{for $u,v\subset\IS$, $u\cap v=\varnothing$, $\ls|u|=\ls|v|=2$ or~$1$},
\]
respectively, and we have
\[
\be{o}\cdot(\vv{u}-\vv{v})=\ls|u\cap o|-\ls|v\cap o|\quad
 \text{whenever $\ls|u|=\ls|v|$}.
\label{eq.A.prod.mixed}
\]

Likewise, the root lattice~$\bD_n$
is the maximal even sublattice of~$\bH_n$:
\[*
\textstyle
\bD_n=\bigl\{\sum_i\Ga_i e_i\in\bH_{n}\bigm|\sum_i\Ga_i=0\bmod2\bigr\}.
\label{eq.Dn}
\]
The discriminant $\discr\bD_n$ is $\Z/2\oplus\Z/2$ (if $n$ is even) or $\Z/4$
(if $n$ is odd); the shortest vectors are
(in the notation $\cl1,\cl2,\cl3$ of~\cite{Conway.Sloane}
for the discriminant classes,
which we underline to avoid confusion)
\[*
\cl2\ni e_i\ \text{for $i\in\IS$}\quad\text{and}\quad
 \cl1,\cl3\ni\be{o}:=\frac12(\vv{o}-\vv{\bar o}),\quad\be{o}^2=\frac{n}4,\quad
 \text{for $o\subset\IS$}
\]
(the class $\be{o}\bmod\bD_n=\cl1$ or~$\cl3$ depends on the parity of~$\ls|o|$),
and we have
%an almost
{\2a}
literal analogue of~\eqref{eq.A.2}:
\[
\be{r}^2-\be{r}\cdot\be{s}=\frac12\ls|r\sdif s|.
\label{eq.D.prod}
\]
If $\2n\ge5$,
the $\RG(\bD_n)$-orbits of integral vectors of square~$m=4$ or~$2$ are
\[*
\orbit4^2:=\bigl\{\pm2e_i\bigr\},\quad
\orbit4:=\bigl\{{\textstyle\sum_{i\in u}\pm e_i}\},\quad
\orbit2:=\bigl\{{\textstyle\sum_{i\in u}\pm e_i}\}\quad
\text{for $u\subset\IS$, $\ls|u|=m$}.
\]
If $\2n=4$, then $\orbit4$ splits into two orbits according to the parity of
the number of the $+$-signs. On the other hand, $\orbit4$ and
$\orbit4^2$ constitute a single $\OG(\bD_4)$-orbit.

We intentionally use the same notation $\be{{}\cdot{}}$ for the shortest
discriminant vectors
in both~$\bA_n$ and~$\bD_n$: in addition to \eqref{eq.A.2},
\eqref{eq.D.prod} \vs. \eqref{eq.H.2}, these vectors have the following
common property:
\[
\sum_{o\subset\IS}\Ga_o\be{o}=\sum_{o\subset\IS}{\2\Ga_o}\vv{o}
 \quad\text{provided that}\quad
 \sum_{o\subset\IS}\Ga_o=0
\label{eq.combination}
\]
for a collection of coefficients $\Ga_o\in\Q$, $o\subset\IS$.
This identity follows from the fact that the coefficients of $e_i$,
$i\in\bar{o}$, differ from those of $e_i$, $i\in o$, by $-1$.
In fact, \eqref{eq.combination} explains also~\eqref{eq.H.2}, \eqref{eq.A.2},
and~\eqref{eq.D.prod}: given two vectors $a,b\in\bH_n\otimes\Q$
of the same length,
$a^2=b^2$, one obviously has $(a-b)^2=2a^2-2(a\cdot b)$.

\subsection{Counts and bounds \via\ blocks}\label{s.blocks}
Assume that an $\RG_h(\bN)$-invariant, \cf.~\eqref{eq.action},
orthogonal decomposition of the form
\[*
\bN\subset B_1\oplus B_2
%\label{eq.blocks}
\]
has been fixed, where $B_r$, $r=1,2$, are certain
positive definite inner product spaces over~$\Q$,
called \emph{blocks}.
%, such that each intersection
%$\bN\cap B_r$ is preserved by $\RG_h(\bN)$.
(Often, we  take for $B_r$ sums of some of
$D_k\otimes\Q$ in~\eqref{eq.Niemeier}, but this is not assumed; in fact, we do not
even assume that $\rank\bN=\dim B_1+\dim B_2$, \cf. \autoref{s.A} below.)
The orthogonal projection $\bN\to B_r$ is denoted by $|_r$, $r=1,2$,
\cf. the similar notation in \autoref{s.Niemeier}.

In view of~\eqref{eq.Rh}, the $\RG_h(\bN)$-invariance of the block
decomposition is equivalent to the assertion that each root
$e\in h^\perp\subset\bN$ is either in~$B_1$ or in~$B_2$; hence, in the
obvious notation, $\RG_h(\bN)=\RG_h(B_1)\times\RG_h(B_2)$.
It follows that
each combinatorial orbit $\orb$ has the property that
\roster
\item\label{i.orb.split}
$\orb=\orb|_1\times\orb|_2=\bigl\{l_1+l_2\bigm|l_r\in\orb|_r,\ r=1,2\bigr\}$ and
\item\label{i.orbit.trans}
%for $\{r,s\}=\{1,2\}$,
the pointwise stabilizer of $\orb|_r$ in $\RG_h(\bN)$ is transitive on $\orb|_s$,
$\{r,s\}=\{1,2\}$.
\endroster
Thus, we can define the \emph{partial count} and \emph{bound} of $\orb$ \via
\[
\cnt_r(\orb):=\bigl|\orb|_r\bigr|,\qquad
\bnd_r(\orb):=\max\ls|\fC_r|,
\label{eq.rel.cb}
\]
the maximum running over all subsets
$\fC_r\subset\orb|_r$ such that, for some
(equivalently, any) fixed
$l_s\in\orb|_s$,
$\{r,s\}=\{1,2\}$,
the set $\bigl\{l_s+l_r\bigm|l_r\in\fC_r\bigr\}$ is admissible.
%the primitive hull in~$\bN$ of the lattice spanned by $h$
%and all $l_r+l_s$, $l_r\in\fC_r$, is root free.
In view of~\iref{i.orbit.trans} above, this latter condition is
indeed independent of the choice of $l_s\subset B_s$.
As an immediate consequence of these definitions and~\iref{i.orb.split},
\iref{i.orbit.trans} above, we have
\[*
\cnt(\orb)=\cnt_1\cnt_2,\qquad
\bnd(\orb)\le\min\{\bnd_1\cnt_2,\cnt_1\bnd_2\}
 =\cnt(\orb)\min_r\frac{\bnd_r}{\cnt_r}.
\]
By induction, we easily extend these relations to any number
of pairwise orthogonal $\RG_h(\bN)$-invariant blocks $B_r$, $r\in I$, arriving at
\[
\cnt(\orb)=\prod_r\cnt_r(\orb),\qquad
\bnd(\orb)\le\cnt(\orb)\min_r\frac{\bnd_r(\orb)}{\cnt_r(\orb)}.
\label{eq.cb}
\]

Since the block decomposition is $\RG_h(\bN)$-invariant, the
following quantities are constant for each $r=1,2$ and each combinatorial
orbit~$\orb$ (where $l\in\orb|_r$):
\[*
\orb|_r^2:=l^2,
%\ l\in\orb|_r,
\qquad
\orb|_r\cdot h|_r:=l\cdot h|_r,
%\ l\in\orb|_r,\quad\text{and}
\qquad
h|_r^2\ (\text{independent of $\orb$}).
\]
Since reflections act identically on the discriminant, the discriminant class
\[*
l\bmod(N\cap B_r)\in\discr(N\cap B_r),\quad l\in\orb|_r,
\]
is also constant within each combinatorial orbit~$\orb$. If this class is
zero, then either $\orb|_r=\{0\}$, or $\orb|_r$ consists of roots, or
$\orb|_r^2=4$, and then necessarily $\orb|_r\cdot h|_r=2$
as the projection of~$\orb$ to the other block is trivial.

For a subset $\fC_r\subset\orb|_r$, we can define the \emph{relative span}
\[*
\spnO\fC_r:=\bigl\{{\textstyle\sum_l\Ga_ll}\bigm|
 \text{$l\in\fC_r$, $\Ga_l\in\Q$, $\textstyle\sum_l\Ga_l=0$}\bigr\}\subset B_r.
%\quad\text{and}\quad
%\spn^0_h\fC_r:=\spn^0\fC_r\oplus\Q h|_r;
\]
Alternatively, $\spnO\fC_r$ can be defined as the vector space spanned by
all differences $l-l_0$, where $l,l_0\in\fC_r$ and $l_0$ is fixed.
Then, necessary (but typically not sufficient) conditions for the
admissibility of
%$\fC_r$
{\2the set $\bigl\{l_s+l_r\bigm|l_r\in\fC_r\bigr\}$ (with $l_s\in\fC_s$
fixed)}
are (see \autoref{lem.intr} and
\autoref{def.admissible})
\begin{align}
\label{eq.intr}
&\text{$l_1^{2}-l_1\cdot l_2=0$ (iff $l_1=l_2$), $2$, $3$, $4$, or $6$
 for all $l_1,l_2\in\fC_r$},\\
\label{eq.span}
&\text{the lattice $\bN\cap\spnO\fC_r$ is root free}.
\intertext{In the special case where $h|_r^2=4$ (or, equivalently, $h\in B_r$),
there is a somewhat stronger restriction still ``local'' with respect to $B_r$:}
\label{eq.spanh}
&\text{the lattice $\bN\cap(\spnO\fC_r+\Q h|_r)$ is root free}.
\end{align}
Clearly, \eqref{eq.spanh} implies~\eqref{eq.span} which, in turn, implies
\emph{part} of~\eqref{eq.intr}: $l_1^{2}-l_1\cdot l_2=1$
if and only if $l_1-l_2$ is a root. The other part of~\eqref{eq.intr},
$l_1^{2}-l_1\cdot l_2\ne5$, {\2can be seen as the ``local''} version of the requirement that
$h-l_1-l_2$ should not be a root.

The fact that conditions~\eqref{eq.intr} and~\eqref{eq.span} are not
sufficient is illustrated by the next lemma, whose proof is not quite local.
The lemma is particularly meaningful if $h|_1=0$ or, more generally,
%$h|_1$ is orthogonal to $\orb|_1$;
$\orb|_1\cdot h|_1=0$;
otherwise, there are better bounds that
are
explored below on a case-by-case basis.

\lemma\label{lem.max.roots}
Consider a block decomposition $\bN\subset B_1\oplus B_2$, let
$\bN_r:=\bN\cap B_r$, and denote by ${\rt_r}\subset\bN_r$ the sublattice generated
by all \emph{integral} roots in $\bN_r$, $r=1,2$.
If, for a combinatorial orbit~$\orb$, we have $\orb|_1\subset{\rt_1}$,
then $\bnd_1(\orb)\le\rank{\rt_1}$\rom; moreover,
the adjacency graph of any admissible set
$\fC_1\subset\orb|_1$ is a Dynkin diagram.
\endlemma

\proof
It follows from the assumption that also $\orb|_2\subset{\rt_2}$.
%, as otherwise the vectors in $\orb|_2$ would be roots in~$\bN$ that are not
%integral.
Let $\fC_1\subset\orb|_1$ be a set as in~\eqref{eq.rel.cb}, so
that, for some fixed $v\in\orb|_2\subset{\rt_2}$, the set
$\fC:=\bigl\{u+v\bigm|u\in\fC_1\bigr\}$ is admissible. By~\eqref{eq.intr}, we have
$u_1\cdot u_2\in\{0,-1,-2\}$ for any distinct $u_1,u_2\in\fC_1$, \ie, the
elements of~$\fC_1$ constitute the vertices of a certain Dynkin diagram,
\latin{a priori} elliptic or affine. However, if
$u_1,\ldots,u_m$ constitute an \emph{affine} Dynkin diagram, there is a
relation $\sum_{i=1}^{m}\Ga_iu_i=0$ with all $\Ga_i>0$, so that
$\Ga:=\sum_{i=1}^{m}\Ga_i\ne0$. It follows that the root
$v=\bigl(\sum_{i=1}^{m}\Ga_i(u_i+v)\bigr)/\Ga$ is in
$\spn\fC$ and $\fC$ is not admissible.

Thus, the set~$\fC_1$ is linearly independent, and the statement
follows.
\endproof

The count $\cnt_1(\orb)$ in \autoref{lem.max.roots} is the total
number of roots in~$\rt_1$, which, for indecomposable root lattices, is as
follows (\cf. $\orbit2$ in \autoref{s.HAD}):
%;
%$\bbC(2,n)=n(n-1)/2$ stand for the binomial coefficients):
\[*
\bA_n\:\ 2\bbC(2,n+1),\quad
\bD_n\:\ 4\bbC(2,n),\quad
\bE_6\:\ 72,\quad
\bE_7\:\ 126,\quad
\bE_8\:\ 240,
\]
where $\bbC(2,n)=n(n-1)/2$ are the binomial coefficients.

\subsection{Combinatorial estimates}\label{s.combinat}
In this section we improve and refine some of the combinatorial bounds
introduced in~\cite{degt:sextics,DIO}. The refinement is based on the extra
condition \eqref{eq.S.span}, which is derived from~\eqref{eq.span} and is
due to the fact that we
only consider primitive sublattices of Niemeier lattices.

Let $\set:=\set_n$ be a finite set, $\ls|\set|=n$.
Given a collection $\fS$ of subsets
of~$\set$, we can consider the $\Q$-vector space
\[*
\spnO\fS:=\bigl\{{\textstyle\sum_s\Ga_s\vv{s}}\bigm|
 \text{$s\in\fS$, $\Ga_s\in\Q$, $\textstyle\sum_s\Ga_s=0$}\bigr\}
 \subset\bH_n\otimes\Q.
\]
The collection~$\fS$
is called \emph{admissible} if
\begin{align}
\label{eq.S.intr}
&\text{$\ls|r\sdif s|\in\Delta:=\{0,4,6,8,12\}$ for each pair $r,s\in\fS$, and}\\
\label{eq.S.span}
&\text{the lattice $\bD_n\cap\spnO\fS$ is root free,
 where $\bD_n\subset\bH_n$ is as in \autoref{s.HAD}}.
\end{align}
As explained in \autoref{s.blocks}, \eqref{eq.S.span} implies \emph{part}
of~\eqref{eq.S.intr}, \viz. the fact that $\ls|r\sdif s|\ne2$. Note though
that our estimates almost do not use the other part, $\ls|r\sdif s|\ne8$.

\remark\label{rem.AD.S}
In view of~\eqref{eq.combination}, the vectors $\vv{s}$ in the definition of
$\spnO\fS$ can be replaced with $\be{s}\in\bA_{n-1}\dual$ or $\be{s}\in\bD_n\dual$,
see \autoref{s.HAD}. This interpretation explains the relevance of admissible
collections: they represent conditions~\eqref{eq.intr}, \eqref{eq.span} that
are necessary for the admissibility of the restriction of a combinatorial
orbit to a block of type $\bA_{n-1}\otimes\Q$ or $\bD_n\otimes\Q$.
\endremark

%\remark\label{rem.invariant}
The interpretation \via\ shortest discriminant vectors $\be{s}\in\bD_n\dual$
implies also that the admissibility
property is invariant under
\[\label{eq.trans}
\aligned
\subeq\label{eq.trans.Sn}
&\text{the natural action of the permutation group $\SG{n}$
on~$\set_n$,}\\
\subeq\label{eq.trans.sdif}
&\text{transformations
$\fS\mapsto\fS\sdif o:=\bigl\{s\sdif o\bigm|s\in\fS\bigr\}$
with $o\subset\set$ fixed},\\
\subeq\label{eq.trans.complement}
&\text{the setwise complement
$\fS\mapsto\fS^-:=\{\set\sminus s\,|\,s\in\fS\}=\fS\sdif\set$}.
\endaligned
\]
Most statements are obvious.
The invariance of~\eqref{eq.S.intr} under \eqref{eq.trans.sdif}
follows from the relation
$(r\sdif o)\sdif(s\sdif o)=r\sdif s$,
and for~\eqref{eq.S.span} we observe that $\fS\mapsto\fS\sdif\{i\}$ is
induced by the reflection~$r_i$ (see \autoref{s.HAD}) which preserves
$\bD_n\subset\bH_n$.

For $0\le m\le n$, define
\[*
\bbA(m,n):=\max\ls|\fS_m|,\qquad \bbD(n):=\max\ls|\fS_*|,
\]
where $\fS_m$ runs over all admissible collections of \emph{$m$-element}
subsets of~$\set$ and $\fS_*$ runs over all admissible collections.
Then, for $1\le m\le n$, we have
\setcounter{subeq}0
\[
\label{eq.bA}
\gathered
\bbA(0,n)=\bbA(1,n)=1,\quad
\bbA(3,6)=\prim3,\quad
\bbA(3,7)=\prim5,\quad
\bbA(3,8)=\prim6,\quad
\bbA(4,8)=\prim8,
\\
\subeq\label{eq.bA.1}
\bbA(m,n)=\bbA(n-m,\2n),\qquad
\subeq\label{eq.bA.2}
\bbA(m,n)\le\bbA(m-1,n-1)+\bbA(m,n-1),\\
\subeq\label{eq.bA.3}
\bbA(m,n)\le\Bigl\lfloor\frac{n}{m}\bbA(m-1,n-1)\Bigr\rfloor,
\endgathered
\]\[
\label{eq.bD}
\gathered
\bbD(0)=\ldots=\bbD(3)=1,\quad
\bbD(4)=\bbD(5)=2,\quad
\bbD(6)=\prim3,\quad
\bbD(7)=\prim5,\quad
\bbD(8)=\prim8,
\\
%\bbD(n)\le\max_{m\ge0}\sum_{s=0}^{6}\bbA(m+2s,n),
\subeq\label{eq.bD.1}
\bbD(n)\le\sum_{m\in\Delta}\bbA(m,n),\qquad
\subeq\label{eq.bD.2}
\bbD(n)\le2\bbD(n-1),
\endgathered
\]
where marked with a $\prim{}$ are the values depending on~\eqref{eq.S.span}.
These assertions are proved below.

\remark\label{rem.AD}
Relations~\eqref{eq.bA},~\eqref{eq.bD}
are used to estimate $\CA_*$, $\CD_*$ inductively,
selecting \emph{at each step} the best bound found. For example, once
$\bbA(2,6)=\bbA(3,6)=3$ has been established, \eqref{eq.bA.3} gives us
$\bbA(3,7)\le\lfloor(7/3)\cdot3\rfloor=7$.
However, using also \eqref{eq.bA.1}, we get the sharp bound
$\bbA(3,7)=\bbA(4,7)\le\lfloor(7/4)\cdot3\rfloor=5$.
\endremark

To make the notation uniform, introduce also
\[*
\bbC(m,n):=\binom{n}{m}\quad\text{and}\quad
\bbP(0):=1,\quad\bbP(n):=2^{n-1}\quad\text{for $n>0$}.
\]
(The count $\bbP(n)$ is interpreted as the number of subsets $s\in\set$ of
the same parity.)
We extend $\bbA(m,n)=\bbC(m,n)=0$ unless $0\le m\le n$ and
$\bbD(n)=\bbP(n)=0$ unless $n\ge0$.

\proof[Proof of~\eqref{eq.bA} and~\eqref{eq.bD}]
Relation~\eqref{eq.bA.1} is the invariance of the
admissibility under \eqref{eq.trans.complement}.
The proofs
of~\eqref{eq.bA.2} and \eqref{eq.bD.2} mimic the
standard proofs of the similar combinatorial identities
for~$\CC_*$ and $\bbP(*)$: we pick a point ${*}\in\set$,
break an admissible set~$\fS$ into the subsets
\[*
\fS-{*}:=\bigl\{s\sminus{*}=s\sdif\{{*}\}\bigm|{*}\in s\in\fS\bigr\}
\quad\text{and}\quad
\fS|_*:=\bigl\{s\bigm|{*}\notin s\in\fS\bigr\},
\]
and estimate the two separately, without any further analysis.

For~\eqref{eq.bA.3}, using
$\ls|\fS-{*}|\le\bbA(m-1,n-1)$ again, we conclude that each point~$*$ is contained in at most
$\bbA(m-1,n-1)$ sets $s\in\fS$.
Counting the total cardinality of the sets $s\in\fS$, we get
$m\ls|\fS|\le n\bbA(m-1,n-1)$, as stated.
For~\eqref{eq.bD.1} we observe that, replacing $\fS$ with
$\fS\sdif o$ for a fixed $o\in\fS$,
see~\eqref{eq.trans.sdif},
we can assume that $\varnothing\in\fS$ and, by~\eqref{eq.S.intr},
$\ls|s|\in\Delta$ for each
$s\in\fS$, upon which we use $\bbA(m,n)$ to
estimate the number of sets of each cardinality separately.

It remains to compute the not immediately obvious
%(those marked with a~$\prim{}$)
values for small values of $m$, $n$.
Below, all uniqueness
assertions are modulo the applicable equivalence given by~\eqref{eq.trans},
\ie, up to~\eqref{eq.trans.Sn} and, for
$\CD_*$-type bounds, also up to~\eqref{eq.trans.sdif}, \eqref{eq.trans.complement}.
We use the term \emph{$(m,n)$-collection} for a collection of
$m$-element subsets of~$\set_n$.

It is easily seen that there are but two maximal (with respect to inclusion)
$(3,6)$-collections satisfying~\eqref{eq.S.intr},
\viz. a pair of disjoint
triplets and
%the one shown in~
\eqref{eq.A.6}:
\[
\label{eq.sets.A.6}
{\subeq\label{eq.A.6}
\miniset
+***...\\
+*..**.\\
-.*.*.*\\
-..!.!!\\
\endminiset
\qquad
\subeq\label{eq.A.7-}
\miniset
+***....\\
-*..**..\\
&*....**\\
&.*.*.*.\\
+.*..*.*\\
-..**..*\\
\endminiset
\qquad
\subeq\label{eq.A.7}
\miniset
* * * . . . . \\
* . . * * . . \\
* . . . . * * \\
. * . * . * . \\
. . * . * . * \\
\endminiset}
\]
In the latter case, the combination shown in the figure is twice a
root; hence, (any) one of the sets (shown as a ``ghost'' in the figure) must
be removed and we arrive at $\bbA(3,6)=3$, realized by a unique
collection.
Then we mimic the proof of \eqref{eq.bA.2}, trying to
construct a larger collection from a pair of smaller ones. The only pair of
extremal
$(2,6)$- and $(3,6)$-collections satisfying~\eqref{eq.S.intr} is
\eqref{eq.A.7-}: it has several prohibited combinations like
the one shown in the figure, and \emph{all} of them cannot be destroyed by
removing just one set. Hence, we need to use the other maximal
$(3,6)$-collection, arriving at \eqref{eq.A.7}.
The passage to~$\bbA(3,8)$ and $\bbA(4,8)$ is similar, but tedious;
we leave details to the reader (or, rather, computer).
%\[
%\subeq\miniset
%+****..\\
%+..****\\
%-!!..!!\\
%-......\\
%\endminiset
%\]
There are but two extremal $(3,8)$-collections, see
\eqref{eq.A.8.3-1}, \eqref{eq.A.8.3-2}, and a unique extremal
$(4,8)$-collection \eqref{eq.A.8.4}.
\[
\label{eq.A.8}
{\subeq\label{eq.A.8.3-1}
\miniset
* * * . . . . . \\
* . . * * . . . \\
* . . . . * * . \\
. * . * . * . . \\
. * . . * . . * \\
. . * . * . * . \\
\endminiset
\quad
\subeq\label{eq.A.8.3-2}
\miniset
* * * . . . . . \\
* . . * * . . . \\
* . . . . * * . \\
. * . * . * . . \\
. * . . * . . * \\
. . * . . . * * \\
\endminiset
\quad
\subeq\label{eq.A.8.4}
\miniset
* * * * . . . . \\
* * . . * * . . \\
* * . . . . * * \\
* . * . * . * . \\
. * . * . * . * \\
. . * * * * . . \\
. . * * . . * * \\
. . . . * * * * \\
\endminiset}
\]

The proof of the $\CD$-type bounds is similar. Up to the new
equivalence,
{\2\viz. \eqref{eq.trans} rather than just the $\SG{n}$-action,}
\eqref{eq.A.6} and
\eqref{eq.A.7} are still the only sets representing
$\bbD(6)=3$ and $\bbD(7)=5$, respectively. The bound $\bbD(8)=8$ is
represented by several classes, including \eqref{eq.A.8.4}.
\endproof

For yet another combinatorial bound we replace the pair $N\ni h$ with
\[*
N:=\bD_{n+1}\subset\bH_{n+1}\quad\text{and}\quad h:=2e_{n+1}
\]
and consider the combinatorial orbit $\orb\subset\orbit4$. Since
$\orb\cdot h=2$, we have
\[*
\orb=\bigl\{\pm e_i\pm e_j\pm e_k+e_{n+1}\bigm|1\le i<j<k\le n\},
\]
so that we can merely speak about square~$3$ vectors
$\pm e_i\pm e_j\pm e_k\in H_n$. (If $n=3$, the set $\orb$ splits into two
combinatorial orbits, but we combine them together.)

Denoting by $\bbT(n):=\max\ls|\fS|$
the maximal cardinality of a set $\fS\subset\orb$ that is admissible
(as in \autoref{def.admissible}) and self-dual, $\fS^*=\fS$,
we have, for $n\ge2$,
\[
\gathered
\bbT(0)=\bbT(1)=\bbT(2)=0,\quad
\bbT(3)=2,\quad
\bbT(4)=4,\quad
\bbT(5)=8,\quad
\bbT(6)=12,\quad
\bbT(7)=18,
\\
\subeq\label{eq.bT.1}
\bbT(n)\le2\Bigl\lfloor\frac13n(n-2)\Bigr\rfloor,\qquad
\subeq\label{eq.bT.2}
\bbT(n)\le2(n-2)+\bbT(n-1).
\endgathered
\label{eq.bT}
\]

\proof[Proof of~\eqref{eq.bT}]
The proof of \eqref{eq.bT.1} follows that of \eqref{eq.bA.3}: we estimate the
total cardinality of the
elements of~$\orb$ regarded as $3$-element subsets of $\{\pm\}\times\IS_n$.
We can take~$2$ out of $\lfloor{\cdot}\rfloor$ since $\bbT(n)$ is obviously
even.
The other bound \eqref{eq.bT.2} is proved similar to
\eqref{eq.bA.2}: we break an admissible set into \emph{three} subsets,
according to the coefficient of~$e_n$, which can be~$0$ or $\pm1$.

In both cases, the r\^{o}le of $\bbA(m-1,n-1)$ should be played by the
bound $b=n-1$ given by
\autoref{lem.max.roots}. However, if roots constituting a Dynkin diagram~$D$
span $\bD_n\otimes\Q$, each component of~$D$ must be of type~$\bD_k$
(including $\bD_3:=\bA_3\subset\bH_3$ or $\bD_2:=2\bA_1\subset\bH_2$);
then, there is a pair $u_1,u_2$ of roots such that $u_1-u_2=2e_i$ and
$\frac12(h+u_1-u_2)\in\spn\fC$ is a root. Thus, we can reduce $b$ down to
$n-2$.

The values $\bbT(n)$ for small values of~$n$ are found by brute force.
%We regard members $s\in\fS$ as $3$-element subsets of $\set:=\{\pm\}\times\IS$.
%Pick a point ${*}=\pm e_i\in\set$ and observe that the set $\fS-{*}$ in~\eqref{eq.S-*}
%can be regarded as an admissible set of roots in the maximal root lattice in
%$\bH_{n-1}\cong(\pm e_i)^\perp\subset\bH_n$ (which is $\bD_{n-1}$, see
%\autoref{s.D} below). Hence, by \autoref{lem.max.roots},
%we have $\ls|\fS-{*}|\le n-1$. Proceeding as in the proof of~\eqref{eq.bA}
%and counting the total cardinality of the sets $s\in\fS$, we arrive at
%$3\ls|\fS|\le\ls|\set|(n-1)=2n(n-1)$.
\endproof

\subsection{Blocks of type $\bA\sb n$}\label{s.A}
In this and next sections, we use the material of \autoref{s.combinat} to
estimate the partial bounds $\bnd_k(\orb)$ for a block~$D_k$ of type~$\bA_n$
or~$\bD_n$. For the three exceptional blocks of type
$\bE_6$, $\bE_7$, $\bE_8$,
the bounds are computed by brute force on the fly, and we do not state the
results here. In fact, we do the same for $\bA_n$, $n\le7$, and
$\bD_n$, $n\le8$, obtaining slightly better bounds in a few cases.

Consider a block $D_k\otimes\Q$ in decomposition~\eqref{eq.Niemeier}
with $D_k$ of type~$\bA_n$.
We assume that $\hk=\sum\eta_ie_i$, $i\in\IS$ and, for each $\eta\in\{\eta_i\}$,
let
\[*
\IS(\eta):=\bigl\{i\bigm|\eta_i=\eta\bigr\}\subset\IS,\quad
\bH(\eta):=\Z\IS(\eta),\quad
B(\eta):=\bH(\eta)\otimes\Q.
\]
The decomposition $D_k\otimes\Q\subset\bigoplus B(\eta)$ is used to
find the count $\cnt_k(\orb)$ and estimate the bound
$\bnd_k(\orb)$, both \via~\eqref{eq.cb},
of a combinatorial orbit~$\orb$.
We reserve the notation $l\in\orbk$ for a ``typical'' element and denote by
$l(\eta)$ its projection to $B(\eta)$.

We use repeatedly the observation that, if either $\hk$ or~$l$ is of
square~$4$, then necessarily $\hk\cdot l=2$, imposing certain restrictions
to the other vector.

For references, we format the rest of this section as a sequence of
numbered claims. Proofs are omitted as all necessary explanations are given
in the text. In a few cases, the ``block-by-block'' bound can be improved
by ruling out certain linear combinations;
this improvement is shown \via\ $(\text{old bound})\mapsto(\text{new bound})^*$.

\claim\label{A.m}
Assume that $\orbk\subset\cl m\ne\cl0$, so that $l=\be{o}$, $\ls|o|=m$.
\roster
\item\label{A.m.u}
If $\hk=\be{u}$, $u\subset\IS$,
then, by~\eqref{eq.A.1}, the two constants
$i:=\ls|u\cap o|$ and $\hk\cdot l$ determine each other and
the restriction $\orbk$. Applying~\eqref{eq.cb} to the partition
$\{B(\eta)\}$, in view of~\eqref{eq.intr}
and~\eqref{eq.A.2} we get
\[*
\cnt_k(\orb)=\bbC(i,\ls|u|)\bbC(m-i,\ls|\bar u|),\quad
\bnd_k(\orb)\le\min
 \bigl\{\bbA(i,\ls|u|)\bbC(m-i,\ls|\bar u|),\bbC(i,\ls|u|)\bbA(m-i,\ls|\bar u|)\bigr\}.
\]
\item\label{A.m.4}
If $h=\hk=\vv{u}-\vv{v}\in\orbit4$,
then $u\subset o$ and $v\subset\bo$,
\cf.~\eqref{eq.A.prod.mixed},
and we have
\[*
\cnt_k(\orb)=\bbC(m-2,n-3),\quad\bnd_k(\orb)\le\bbA(m-2,n-3).
\]
\item\label{A.m.2}
If $\hk=\vv{u}-\vv{v}\in\orbit2$ is a root, then
\[*
\alignedat3
\cnt_k(\orb)&=\bbC(m-2,n-1),&\quad\bnd_k(\orb)&\le\bbA(m-2,n-1),&
 \quad&\text{if $\hk\cdot l=0$ and $u,v\subset o$},\\
\cnt_k(\orb)&=\bbC(m,n-1),&\quad\bnd_k(\orb)&\le\bbA(m,n-1),&
 \quad&\text{if $\hk\cdot l=0$ and $u,v\subset\bo$},\\
\cnt_k(\orb)&=\bbC(m-1,n-1),&\quad\bnd_k(\orb)&\le\bbA(m-1,n-1),&
 \quad&\text{if $\hk\cdot l=\pm1$}.
\endalignedat
\]
\endroster
\endclaim

\claim\label{A.4}
Assume that $\orbk\subset\orbit4$, so that
$l:=\vv{r}-\vv{s}$, $\ls|r|=\ls|s|=2$.
\roster
\item\label{A.4.u}
If $\hk=\be{u}$, $u\subset\IS$, then $r\subset u$,
$s\subset\bar u$, see~\eqref{eq.A.prod.mixed},
and, by~\eqref{eq.intr} and~\eqref{eq.H.2},
\[*
\cnt_k(\orb)=\bbC(2,\ls|u|)\bbC(2,\ls|\bar u|),\quad
\bnd_k(\orb)\le\min
 \bigl\{\bbA(2,\ls|u|)\bbC(2,\ls|\bar u|),\bbC(2,\ls|u|)\bbA(2,\ls|\bar u|)\bigr\}.
\]
\item\label{A.4.4}
If $h=\hk=\vv{u}-\vv{v}\in\orbit4$, then
\[*
\alignedat3
\cnt_k(\orb)&=\bbC(2,n-3),&\quad\bnd_k(\orb)&\le\bbA(2,n-3),&
 \quad&\text{if $r=u$ or $s=v$},\\
\cnt_k(\orb)&=8\bbC(2,n-3),&\quad\bnd_k(\orb)&\le4(n-4),&
 \quad&\text{if $\ls|s\cap u|=\ls|r\cap v|=1$};
\endalignedat
\]
for the latter, $l(0)$ is a root in
$\bA_{n-4}\subset\bH(0)$
and we use \autoref{lem.max.roots}.
\item\label{A.4.2}
If $\hk=\vv{u}-\vv{v}\in\orbit2$, then $u\subset r$,
$v\subset s$ and, similar to case~\iref{A.4.4},
\[*
\cnt_k(\orb)=2\bbC(2,n-1),\quad\bnd_k(\orb)\le n-2.
\]
\endroster
\endclaim

\claim\label{A.2}
Assume that $\orbk\subset\orbit2$,
so that $l:=\vv{r}-\vv{s}$, $\ls|r|=\ls|s|=1$.
Assume also that $\hk\cdot l\ne0$, as otherwise
\autoref{lem.max.roots} applies to (the components of) $\hk^\perp$.
\roster
\item\label{A.2.u}
If $\hk=\be{u}$, $u\subset\IS$, then $r\subset u$, $s\subset\bar u$
(or \latin{vice versa}), see~\eqref{eq.A.prod.mixed}; hence,
by~\eqref{eq.intr}, \eqref{eq.H.prod}, and~\eqref{eq.bA},
\[*
\cnt_k(\orb)=\ls|u|\ls|\bar u|,\quad
\bnd_k(\orb)\le\min
 \bigl\{\ls|u|,\ls|\bar u|\bigr\}.
\]
\item\label{A.2.4}
If $h=\hk=\vv{u}-\vv{v}\in\orbit4$, then
$r\subset u$, $s\subset v$ and
\[*
\cnt_k(\orb)=4,\quad\bnd_k(\orb)\le2\mapsto\prim1.
\]
Indeed, any pair of distinct vectors in an admissible set $\fC_k\subset\orbk$
is of the form $l$, $h-l$, and then the root $l\in\bN$ is in
$\spnO\fC_k\oplus\Q h$, see \eqref{eq.spanh}.
\item\label{A.2.2}
If $\hk=\vv{u}-\vv{v}\in\orbit2$, then, by~\eqref{eq.bA},
\[*
\alignedat3
\cnt_k(\orb)&=n-1,&\quad\bnd_k(\orb)&=1,&
 \quad&\text{if $\hk\cdot l=\pm1$ (two orbits each)},\\
\cnt_k(\orb)&=1,&\quad\bnd_k(\orb)&=1,&
 \quad&\text{if $\hk\cdot l=\pm2$}.
\endalignedat
\]
\endroster
\endclaim

%\subsection{Conics split into two roots}\label{s.roots}
%
%$1_o,\bold{1}_o,[o]$

\subsection{Blocks of type $\bD\sb n$}\label{s.D}
Consider a block $D_k\otimes\Q$ in decomposition~\eqref{eq.Niemeier}
with $D_k$ of type~$\bD_n$ and
proceed as in \autoref{s.A}, using the notation therein.
Since the group $\OG(\bD_n)$ contains the reflection against any of $e_i$, in
the expression
$\hk=\sum\eta_ie_i$
we can and always do assume that all $\eta_i\ge0$. Then, the description of
combinatorial orbits is still quite ``combinatorial'' in its nature.

\claim\label{D.2d}
Assume that $\orbk\subset\cl2$.
Then,
from~\eqref{eq.intr} and~\eqref{eq.bA}, we have
\[*
\alignedat3
\cnt_k(\orb)&=2\bigl|\IS(0)\bigr|,&\quad\bnd_k(\orb)&\le2,&
 \quad&\text{if $\hk\cdot l=0$},\\
\cnt_k(\orb)&=\bigl|\IS(\eta)\bigr|,&\quad\bnd_k(\orb)&=1,&
 \quad&\text{if $\eta:=\bigl|\hk\cdot l\bigr|>0$}.
\endalignedat
\]
\endclaim

\claim\label{D.1,3}
Assume that $\orbk\subset\cl1$ or~$\cl3$, so that $l=\be{o}$, $o\subset\IS$.
Then $\hk\notin\orbit4^2$.
\roster
\item\label{D.1,3.1u}
If $\hk=0$ or $\hk\in\cl2$ or $\hk\in\orbit4$, \ie,
$\hk=\vv{u}$, $\ls|u|\in\{0,1,4\}$, then
\[*
\cnt_k(\orb)=\bbP(n-\ls|u|),\quad\bnd_k(\orb)\le\bbD(n-\ls|u|).
\]
\item\label{D.1,3.1,3}
If $\hk=\be{\IS}\in\cl1$ or~$\cl3$, then $\ls|o|$ is determined by
$\hk\cdot l=\const$ and
\[*
\cnt_k(\orb)=\bbC(\ls|o|,n),\quad\bnd_k(\orb)\le\bbA(\ls|o|,n).
\]
\item\label{D.1,3.2}
If $\hk=\vv{u}\in\orbit2$, then
\[*
\alignedat3
\cnt_k(\orb)&=2\bbP(n-2),&\quad\bnd_k(\orb)&\le2\bbD(n-2),&
 \quad&\text{if $\hk\cdot l=0$},\\
\cnt_k(\orb)&=\bbP(n-2),&\quad\bnd_k(\orb)&\le\bbD(n-2),&
 \quad&\text{if $\hk\cdot l=\pm1$}.
\endalignedat
\]
\endroster
\endclaim

\claim\label{D.4-2}
Assume that $\orbk\subset\orbit4^2$.
Then $\hk=\vv{u}$, $\ls|u|\in\{1,2,4\}$, and
\[*
\cnt_k(\orb)=\ls|u|,\quad\bnd_k(\orb)\le\ls|u|\mapsto\prim1.%\ls|u|.
\]
Indeed, if $l_1\ne l_2$, then $\2\frac12(l_1-l_2)$ would be a root.
\endclaim

\claim\label{D.4}
Assume that $\orbk\subset\orbit4$, so that
$l=\sum_{i\in o}\pm e_i$, $\ls|o|=4$.
Then $\hk\notin\cl0,\cl2$.
\roster
\item\label{D.4.1,3}
If $\hk=\be{\IS}\in\cl1$ or~$\cl3$, then, by~\eqref{eq.intr}
and~\eqref{eq.D.prod},
\[*
\cnt_k(\orb)=\bbC(4,n),\quad\bnd_k(\orb)\le\bbA(4,n).
\]
\item\label{D.4.4-2}
If $\hk=2e_i\in\orbit4^2$, then, by the definition of $\bbT(*)$,
\[*
\cnt_k(\orb)=8\bbC(3,n-1),\quad
%\bnd_k(\orb)\le\bbT(n-1).
{\2\bnd_k(\orb)=\bbT(n-1)}.
\]
\item\label{D.4.4}
If $\hk=\vv{u}\in\orbit4$, $\ls|u|=4$, then there are two orbits $\orbk$:
\[*
\alignedat3
\cnt_k(\orb)&=4,&\quad\bnd_k(\orb)&\le4\mapsto\prim1,&
 \quad&\text{if $o=u$},\\
\cnt_k(\orb)&=4\bbC(2,4)\bbC(2,n-4),&\quad\bnd_k(\orb)&\le\bbC(2,4)(n-4),&
 \quad&\text{if $\ls|o\cap u|=2$}.
\endalignedat
\]
In the former case, the sign pattern of each vector $l\in\orbk=\orb$ is ${+}{+}{+}{-}$
(up to permutation); hence, $\2\frac12(l_1-l_2)$ would be a root whenever
$l_1\ne l_2$.
In the latter case, \autoref{lem.max.roots} is applied to $\bH(0)$.
\item\label{D.4.2}
If $\hk=\vv{u}\in\orbit2$, $\ls|u|=2$, then $u\subset o$ and, by
\autoref{lem.max.roots},
\[*
\cnt_k(\orb)=4\bbC(2,n-2),\quad\bnd_k(\orb)\le n-2.
\]
\endroster
\endclaim

\claim\label{D.2}
Assume that $\orbk\subset\orbit2$, so that
$l=\sum_{i\in o}\pm e_i$, $\ls|o|=2$. Assume also that $\hk\cdot l\ne0$,
as otherwise \autoref{lem.max.roots} applies to (the components of)
$\hk^\perp$.
\roster
\item\label{D.2.2d}
If $\hk=e_i\in\cl2$,
then $i\in o$ and, by~\eqref{eq.intr} and~\eqref{eq.D.prod},
\[*
\cnt_k(\orb)=2(n-1),\quad\bnd_k(\orb)\le2.
\]
\item\label{D.2.1,3}
If $\hk=\be{\IS}\in\cl1$ or~$\cl3$, then $\hk\cdot l=\pm{\21}=\const$ and,
by~\eqref{eq.intr},
\[*
\cnt_k(\orb)=\bbC(2,n),\quad\bnd_k(\orb)\le\bbA(2,n).
\]
\item\label{D.2.4-2}
if $\hk=2e_i\in\orbit4^2$,
then $i\in o$ and,
similar to \autoref{A.2}\iref{A.2.4},
%by~\eqref{eq.intr} and~\eqref{eq.A.prod},
\[*
\cnt_k(\orb)=2(n-1),\quad\bnd_k(\orb)\le2\mapsto\prim1.
\]
\item\label{D.2.4}
If $\hk=\vv{u}\in\orbit4$, $\ls|u|=4$, then $o\subset u$ and,
similar to \autoref{A.2}\iref{A.2.4},
%by~\eqref{eq.intr},
\[*
\cnt_k(\orb)=\bbC(2,4),\quad\bnd_k(\orb)\le\bbA(2,4)=2\mapsto\prim1.
\]
\item\label{D.2.2}
If $\hk=\vv{u}\in\orbit2$, $\ls|u|=2$, then there are four orbits:
\[*
\alignedat3
\cnt_k(\orb)&=4(n-2),&\quad\bnd_k(\orb)&\le4\mapsto\prim3,&
 \quad&\text{if $\hk\cdot l=\pm1$},\\
\cnt_k(\orb)&=1,&\quad\bnd_k(\orb)&=1,&
 \quad&\text{if $\hk\cdot l=\pm2$}.
\endalignedat
\]
In the formar case, letting $u=\{1,2\}$, a $4$-element admissible set must
be of the form $l_{1,2}=e_1\pm e_3$, $l_{3,4}=e_2\pm e_4$, and
$\frac12(l_1+l_2-l_3-l_4)$ is a root.
\endroster
\endclaim

\subsection{Special $\bD$-type blocks}\label{s.D.special}
Consider a block decomposition $\bN\subset B_1\oplus B_2$ with
$B_1=D_k\otimes\Q$, see \eqref{eq.Niemeier}, and $D_k$ of type~$\bD_n$.
This $\bD$-type block~$B_1$ is called \emph{special} (with respect to a
combinatorial orbit~$\orb$) if either
\roster*
\item
$\orb|_1\subset\cl2$ and $\orb|_1\cdot h|_1=0$, see \autoref{D.2d}, or
\item
$\orb|_1\subset\orbit2$, $h|_1\in\cl2$, and $\orb|_1\cdot h|_1\ne0$,
see \autoref{D.2}\iref{D.2.2d}.
\endroster
In both cases, $\bnd_1(\orb)\le2$, so that \eqref{eq.cb} would typically
result in the rough estimate $\bnd(\orb)\le2\cnt_2(\orb)$. Often, this can be
improved by means of the following lemma.

\lemma\label{lem.D.special}
If the block~$B_1$ in a block decomposition $\bN\subset B_1\oplus B_2$ is
special with respect to a
combinatorial orbit~$\orb$, then $\bnd(\orb)\le\cnt_2(\orb)+\bnd_2(\orb)$.
\endlemma

\proof
Let $\fC\subset\orb$ be an admissible set. For each $l_2\in\orb|_2$,
there are at most two vectors $l_1+l_2\in\fC$ and, if there are two, they are
of the form $\pm e_i+{\const}$, so that their difference is $2e_i$ for some
generator $e_i\in B_1$. If two distinct vectors $2e_i\ne2e_j$ could be obtained
in this way, then their half-sum would be a root in $\spn\fC$.
Hence, at most one vector $2e_i$ can be obtained, and this can be done
starting from at most $\bnd_2(\orb)$ vectors $l_2\in\orb|_2$,
so that we arrive
at
\[*
\bnd(\orb)\le2\cdot\bnd_2(\orb)+1\cdot\bigl(\cnt_2(\orb)-\bnd_2(\orb)\bigr)
 =\cnt_2(\orb)+\bnd_2(\orb).
\qedhere
\]
\endproof

\autoref{lem.D.special} can be iterated: if in a decomposition
$\bN\subset B_0\oplus B_1\oplus\ldots\oplus B_n$ all blocks but~$B_0$
are special, then
\[
\bnd(\orb)\le\bnd_0(\orb)+\sum_{k=1}^{n-1}\prod_{i=0}^{k}\cnt_i(\orb).
\label{eq.D.special}
\]
The best estimate is obtained if the special blocks~$B_k$, $k=1,\ldots,n$,
are ordered by the increasing of $\cnt_k(\orb)$.

\section{Niemeier lattices with many roots}\label{S.many.roots}

The goal of this section is the following statement, ruling out $18$ out of
the $24$ Niemeier lattices, \viz. those with many roots.

\theorem[see \autoref{proof.many}]\label{th.many}
Let $\bN:=\N(D)$ be a Niemeier lattice generated over~$\Q$ by a root
lattice~$D$ \emph{other than}
%the Leech lattice~$\Leech$ or lattice $\N(D)$ with
\[*
6\bD_4,\ 6\bA_4,\ 8\bA_3,\ 12\bA_2,\ \text{\rm or}\ 24\bA_1.
\]
Then, $\ls|\fC|<720$
for any $4$-polarization $h\in\bN$ and any admissible set
$\fC\subset\Orb_h$.
\endtheorem

\subsection{Homogeneous chains of admissible sets}\label{s.chains}
Consider a $4$-polarized Niemeier lattice $\bN\ni h$. For a subset
$\sS\subset\Orb_h$ and subgroup $G\subset\OG_h(\bN)$, denote
\[*
\stab(\sS;G):=\text{the setwise stabilizer of $\sS$ in $G$}.
\]
If $G=\OG_h(\bN)$, the group is omitted from the notation, abbreviating it to
$\stab\sS$.

Fix a pair $\sQ\subset\sS\subset\Orb_h$ of subsets,
non necessarily admissible or saturated,
let $O$ be the orbit of~$\sQ$ under $G:=\stab\sS$,
and assume that the pair is
\emph{homogeneous} in the sense that the \emph{multiplicity}
\[
\label{eq.transitive}
m:=m(l):=\#\bigl\{\sQ'\in O\bigm|l\in\sQ'\bigr\}=\const,\quad l\in\sS,
%\text{the group $G:=\stab\sS$ acts transitively on~$\sS$}.
\]
is constant throughout~$\sS$.
(\latin{E.g.}, this is obviously the case if either $G$ is transitive on $\sS$
or all elements $\sQ'\in O$ are pairwise disjoint and $\bigcup\sQ'=\sS$.)
%Pick another subset $\sQ\subset\sS$ and denote by~$O$ its $G$-orbit.
Denote
\[*
s:=\ls|\sS|,\quad q:=\ls|\sQ|,\quad o:=\ls|O|,\quad\text{so that}\quad
qo=sm.
\]
%where the \emph{multiplicity} $m$ is the number of elements of~$O$
%containing any fixed conic $l\in\sS$; this number is constant due
%to~\eqref{eq.transitive}.
The following obvious lemma on group actions is crucial for our computation.

\lemma\label{lem.chain}
In the notation above, for any subset $\fC\subset\sS$ of cardinality
$n:=\ls|\fC|$, there is an element $\sQ'\subset O$ such that
$\ls|\fC\cap\sQ'|\ge\lceil n\cdot q/s\rceil$.
\done
\endlemma

Recall that our ultimate goal is finding all large admissible/geometric
subsets $\fC\subset\Orb_h$, and this task reduces to the computation of the
sets of the form
\[
\Bnd_n(\sS):=
\bigl\{\text{$\sS$-saturated admissible sets $\fC\subset\sS$}\bigm|
 \ls|\fC|\ge n\bigr\}/\stab\sS.
\label{eq.Bnd}
\]
The latter is computed by brute force, adding conics one-by-one to increase
the rank by one unit at a time. (Certainly, we make use of
the symmetry groups as much as possible, but these
technicalities are discussed in \tabref{s.brute}.) Experimentally, the
number of (orbits of) admissible sets $\fC\subset\sS$ is relatively small
if $\fC$ is either small or close to the maximum, whereas it may grow quite
large in the middle of the range. \autoref{lem.chain} lets us break this
computation into several steps and avoid the above ``middle'':
when computing $\Bnd_n(\sS)$, we can start with
$\Bnd_{\lceil np/s\rceil}(\sQ)$.

\warning\label{wrn.stab}
It may happen that the group $\stab(\sQ;G)$ is smaller than $\stab\sQ$.
In this case, we have two options:
\roster
\item\label{stab.1}
break $(\stab\sQ)$-orbits $[\fC]\in\Bnd_*(\sQ)$
into $\stab(\sQ;G)$-orbits before continuing the computation, or
\item\label{stab.2}
extend $\fC$ not to $\sS$ only, but to all elements of the
$(\stab\fC)$-orbit of~$\sS$.
\endroster
We use both approaches, \cf. \texttt{select} in \tabref{s.code}.
\endwarning

Often, a longer \emph{homogeneous chain}
$\sQ_1\subset\sQ_2\subset\ldots\subset\sQ_N:=\sS$
needs to be used to speed up the computation. Typically,
$\sQ_k$ is the union of several elements of the $(\stab\sQ_{k+1})$-orbit of
$\sQ_{k-1}$.
Details are explained in~\cite{degt:800.tab} case by case.

\subsection{Proof of \autoref{th.many}}\label{proof.many}
We consider Niemeier lattices one-by-one,
using their description in
\cite[Table 16.1]{Conway.Sloane}. For each lattice~$\bN$, we list the
$\OG(\bN)$-orbits of square~$4$ vectors $h\in\bN$ and, for a
representative~$h$ of each orbit, we list orbits and combinatorial orbits
(see \autoref{s.Niemeier}) of conics. Upon this, in most cases it remains to
compute the partial counts and bounds using \autoref{s.A} and \autoref{s.D}
and apply~\eqref{eq.cb} to arrive at $\bnd(\Orb_h)<720$.
We refer to~\tabrefs{\autoref*{S.very.many}, \autoref*{S.many}}
for a full account of this computation and to
\autoref{s.ex.A24}, \autoref{s.ex.2D12} below for several simple examples worked
out in detail.

Occasionally, the bounds need to be reduced by \autoref{lem.D.special}, \cf.
\autoref{s.ex.4D6} below.

In eleven cases, \viz.
\roster*
\item
one polarization for $\N(2\bA_7\oplus2\bD_5)$, see \tabref{s.2A7_2D5-1},
\item
two polarizations for $\N(4\bA_6)$,
see \tabrefs{\autoref*{s.4A6-1}, \autoref*{s.4A6-2}}, and
\item
all eight polarizations for $\N(4\bA_5\oplus\bD_4)$, see \tabref{s.4A5_D4},
\endroster
the exact bounds $\bnd(\orb)$
need to be computed by brute force. This is done step-by-step as explained in
\autoref{s.chains};
the precise choices making the computation reasonably fast are found
in the relevant references above.
\qed

\remark
Probably, \latin{a posteriori} the computation in the last paragraph of the
proof can also be ``explained'' in the spirit of \autoref{lem.D.special}.
We leave these attempts to the reader since, in the subsequent sections, hard
computer aided computation is still unavoidable.
\endremark

\subsection{Example: the lattice $\N(\bA\sb{24})$}\label{s.ex.A24}
The lattice is the index~$5$ extension of~$\bA_{24}$
by the discriminant classes that are multiples of
$\cl5\in\Z/25=\discr\bA_{24}$, see \autoref{s.A}. It is immediate that, up to
$\OG(\bN)$, there are but two square~$4$ vectors $h\in\bN$.

\subsubsection{Case $1$\rom: $h=\vv{u}-\vv{v}\in\orbit4$, $\ls|u|=\ls|v|=2$}
There are five combinatorial orbits, constituting three orbits.
\roster
\item
A pair of dual combinatorial orbits $\orb\subset\cl{20}$, $\orb^*\subset\cl5$:
by \autoref{A.m}\iref{A.m.4} and \eqref{eq.bA}, we have
$\bnd(\orb)=\bnd({\orb^*})\le\bbA(18,21)\le70$.
%\newline
%%\item
%the dual orbit is $\orb^*\subset\cl5$, with the same bound $\bnd(\orb^*)\le70$.
\item
A pair of dual combinatorial orbits
$\orb\subset\orbit4$, with
%a typical vector
$l$ as in the first equation of
\autoref{A.4}\iref{A.4.4}: we have $\bnd(\orb)\le\bbA(2,21)\le10$ for each.
\item
A self-dual combinatorial orbit
$\orb\subset\orbit4$, with
%a typical vector
$l$ as in the second equation of
\autoref{A.4}\iref{A.4.4}: we have $\bnd(\orb)\le4\cdot20=80$.
\endroster
Summarizing, $\bnd(\Orb_h)\le2\cdot70+2\cdot10+80=240$.

\subsubsection{Case $2$\rom: $h=\be{u}\in\cl{20}$, $\ls|u|=20$}
There are but two dual orbits, consisting of a single combinatorial orbit each.
\roster
\item
$\orb\subset\cl{20}$:
by \autoref{A.m}\iref{A.m.u} with $i=18$ and~\eqref{eq.bA}, we have
\[*
{\2\bnd(\orb)}\le\min\{\bbA(18,20)\bbC(2,5),\bbC(18,20)\bbA(2,5)\}
 \le\min\{100,380\}=100.
\]
\item
$\orb\subset\orbit4$: by \autoref{A.4}\iref{A.4.u}, we get the same bound
$\bnd(\orb)\le100$.
\endroster
Summarizing, $\bnd(\Orb_h)\le100+100=200$.

\subsection{Example: the lattice $\N(2\bD\sb{12})$}\label{s.ex.2D12}
The lattice is the index~$4$ extension of $2\bD_{12}$ by the
discriminant classes
$\cl1\oplus\cl2$, $\cl2\oplus\cl1$, and $\cl3\oplus\cl3$,
%$\cl1\oplus\cl2,\cl2\oplus\cl1,\cl3\oplus\cl3\in(\Z/2\oplus\Z/2)^2=\discr2\bD_{12}$,
see \autoref{s.D}.
Up to $\OG(\bN)$, there are four square~$4$ vectors $h\in\bN$.

\subsubsection{Case $1$\rom: $h|_1=2e_i\in\orbit4^2$ and $h|_2=0$}
There are three combinatorial orbits, all self-dual, constituting three
separate orbits.
\roster
\item
$\orb|_1\subset\orbit4$ and $\orb|_2=\{0\}$:
by \autoref{D.4}\iref{D.4.4-2} and~\eqref{eq.bT},
%we have
$\bnd(\orb)\le\bbT(11)\le66$.
\item
$\orb|_1\subset\orbit2$ and $\orb|_2\subset\orbit2$: we have
\[*
\alignedat3
\cnt_1(\orb)&=22,&\quad\bnd_1(\orb)&=1&
 \quad&\text{by \autoref{D.2}\iref{D.2.4-2} and}\\
\cnt_2(\orb)&=264,&\quad\bnd_2(\orb)&\le12&
 \quad&\text{by \autoref{lem.max.roots}},
\endalignedat
\]
resulting in $\bnd(\orb)\le264$.
\item
$\orb|_1\subset\cl2$ and $\orb|_2\subset\cl1$: since
$\cnt_1(\orb)=\bnd_1(\orb)=1$ by \autoref{D.2d},
from~\eqref{eq.bD} and \autoref{D.1,3}\iref{D.1,3.1u} we have
$\bnd(\orb)=\bnd_2(\orb)\le\bbD(12)\le128$.
\endroster
Summarizing, $\bnd(\Orb_h)\le66+264+128=458$.

\subsubsection{Case $2$\rom: $h|_1=\vv{u}\in\orbit2$ and
$h|_2=\vv{v}\in\orbit2$, $\ls|u|=\ls|v|=2$}
There are six orbits.
\roster
\item
$\orb|_1\subset\cl1$ and $\orb|_2\subset\cl2$ (or \latin{vice versa}, two
combinatorial orbits): as
$\cnt_2(\orb)=2$ by \autoref{D.2d}, % yields
$\bnd(\orb)\le2\bnd_1(\orb)\le2\bbD(10)\le64$ by
\autoref{D.1,3}\iref{D.1,3.2}.
\item
$\orbk\subset\orbit2$ and $\orbk\cdot\hk=1$ for $k=1,2$:
we have $\cnt_k(\orb)=40$ and $\bnd_k(\orb)\le3$ by
\autoref{D.2}\iref{D.2.2};
hence, $\bnd(\orb)\le120$ by~\eqref{eq.cb}.
\item\label{i.D.2.3}
$\orb|_1\in\orbit4$ and $\orb|_2=\{0\}$ (or \latin{vice versa}, two
combinatorial orbits): in view of \autoref{D.4}\iref{D.4.2}, we have
$\bnd(\orb)=\bnd_1(\orb)\le10$.
\item
The dual of~\iref{i.D.2.3}, with $\orb|_1\in\orbit2$ and $\orb|_2=\{h|_2\}$
(or \emph{vice versa}).
\item\label{i.D.2.5}
$\orb|_1\in\orbit4^2$ and $\orb|_2=\{0\}$ (or \latin{vice versa}, two
combinatorial orbits): in view of \autoref{D.4-2}, we have
$\bnd(\orb)=\bnd_1(\orb)\le1$.
\item\label{i.D.2.6}
The dual of~\iref{i.D.2.5}, with $\orb|_1\in\orbit2$ and $\orb|_2=\{h|_2\}$
(or \emph{vice versa}).
\endroster
Summarizing,
$\bnd(\Orb_h)\le2\cdot64+120+2\cdot2\cdot10+2\cdot2\cdot1=292$.

\subsubsection{Case $3$\rom: $h|_1=\vv{u}\in\orbit4$, $\ls|u|=4$, and $h|_2=0$}
Each of the five orbits consists of a single combinatorial orbit.
\goodbreak
\roster
\item
$\orb|_1\subset\orbit4$ and $\orb|_2=\{0\}$: we have
$\bnd(\orb)\le\bbC(2,4)\cdot8=48$ by \autoref{D.4}\iref{D.4.4}.
\goodbreak
\item
$\orb|_1\subset\cl1$ and $\orb|_2\subset\cl2$:
%using \autoref{D.1,3}\iref{D.1,3.1u}
%we find that
%$\cnt_1(\orb)=\bbP(8)=128$ and $\bnd_1(\orb)\le\bbD(8)\le8$, so that
%\autoref{lem.D.special} results in $\bnd(\orb)\le128+8=136$.
we have
\[*
\alignedat3
\cnt_1(\orb)&=128,&\quad\bnd_1(\orb)&\le\bbD(8)\le8&
 \quad&\text{by \autoref{D.1,3}\iref{D.1,3.1u} and}\\
\cnt_2(\orb)&=24,&\quad\bnd_2(\orb)&\le2&
 \quad&\text{by \autoref{D.2d}},
\endalignedat
\]
resulting in $\bnd(\orb)\le\min\{128\cdot2,24\cdot8\}=256$.
Since $B_2$ is a special $\bD$-type block, \autoref{lem.D.special} reduces
this down to $\bnd(\orb)\le128+8=136$.
\item
$\orb|_1\subset\orbit2$ and $\orb|_2\subset\orbit2$: we have
\[*
\alignedat3
\cnt_1(\orb)&=6,&\quad\bnd_1(\orb)&\le1&
 \quad&\text{by \autoref{D.2}\iref{D.2.4} and}\\
\cnt_2(\orb)&=264,&\quad\bnd_2(\orb)&\le12&
 \quad&\text{by \autoref{lem.max.roots}},
\endalignedat
\]
resulting in $\bnd(\orb)\le\min\{6\cdot12,264\cdot1\}=72$.
\item\label{i.D.3.4}
$\orb|_1\subset\orbit4^2$ and $\orb|_2=\{0\}$: we have
$\bnd(\orb)=1$ by \autoref{D.4-2}.
\item
the dual of~\iref{i.D.3.4}, with $\orb|_1\subset\orbit4$
as in the first equation of \autoref{D.4}\iref{D.4.4}.
\endroster
Summarizing,
$\bnd(\Orb_h)\le48+256+72+2\cdot1=314$, reducible down to $258$.

\subsubsection{Case $4$\rom: $h|_1=e_i\in\cl2$ and $h|_2=\vv{\IS}\in\cl1$}
There are six orbits, each consisting of a single combinatorial orbit; they
split into three pairs of dual ones.
\roster
\item\label{i.D.4.1}
$\orb|_1=\{2h|_1\}\subset\orbit4^2$ and $\orb|_2=\{0\}$: we have $\bnd(\orb)=1$ by
\autoref{D.4-2}.
\item
the dual of~\iref{i.D.4.1}, with $\orb|_1=\{-h|_1\}\subset\cl2$
and $\orb|_2=\{h|_2\}\subset\cl1$.
\item\label{i.D.4.3}
$\orb|_1\subset\orbit2$ and $\orb|_2\subset\orbit2$:
we have
\[*
\alignedat3
\cnt_1(\orb)&=22,&\quad\bnd_1(\orb)&\le2&
 \quad&\text{by \autoref{D.2}\iref{D.2.2d} and}\\
\cnt_2(\orb)&=66,&\quad\bnd_2(\orb)&\le\bbA(2,12)\le6&
 \quad&\text{by \autoref{D.2}\iref{D.2.1,3}},
\endalignedat
\]
resulting in $\bnd(\orb)\le\min\{22\cdot6,66\cdot2\}=132$.
Since $B_1$ is a special $\bD$-type block, \autoref{lem.D.special} reduces
this down to $\bnd(\orb)\le66+6=72$.
\item
the dual of~\iref{i.D.4.3}, with $\orb|_1\subset\cl2$ and
$\orb|_2\subset\cl1$.
\item\label{i.D.4.5}
$\orb|_1\subset\cl2$ and $\orb|_2\subset\cl1$:
since $\cnt_1(\orb)=\bnd_1(\orb)=1$ by \autoref{D.2d},
using~\eqref{eq.bA} and \autoref{D.1,3}\iref{D.1,3.1,3} we arrive at
$\bnd(\orb)=\bnd_2(\orb)\le\bbA(4,12)\le48$.
\item
the dual of~\iref{i.D.4.5}, with $\orb|_1=\{0\}$ and $\orb|_2\subset\orbit4$.
\endroster
Summarizing,
$\bnd(\Orb_h)\le2\cdot1+2\cdot132+2\cdot48=362$,
reducible down to $242$.

\subsection{Example: the lattice $\N(4\bD\sb6)$}\label{s.ex.4D6}
The lattice is the index~$16$ extension of $4\bD_6$ by the linear
combinations of all even permutations of
$\cl0\oplus\cl1\oplus\cl2\oplus\cl3$. Up to $\OG(\bN)$, there are five
square~$4$ vectors $h\in\bN$. We consider one, $h|_1\in\orbit4^2$ and
$\hk=0$ for $k\ge2$,
\ie, the only class where,
unlike \autoref{s.ex.2D12}, \autoref{lem.D.special} makes a difference.

All combinatorial orbits are self-dual and split into four orbits.
Orbits~\iref{i.D6.2} and~\iref{i.D6.3} below consist of three combinatorial
orbits each, obtained from the one shown by all cyclic permutations
of the indices $2,3,4$.
\roster
\item\label{i.D6.1}
$\orb|_1\subset\orbit4$, $\orbk=\{0\}$ for $k=2,3,4$:
by \autoref{D.4}\iref{D.4.4-2}, $\bnd(\orb)=\bbT(5)=8$.
\item\label{i.D6.2}
$\orb|_1\subset\orbit2$, $\orb|_2\subset\orbit2$, $\orb|_3=\orb|_4=\{0\}$:
we have
\[*
\alignedat3
\cnt_1(\orb)&=10,&\quad\bnd_1(\orb)&=1&
 \quad&\text{by \autoref{D.2}\iref{D.2.4-2} and}\\
\cnt_2(\orb)&=60,&\quad\bnd_2(\orb)&\le6&
 \quad&\text{by \autoref{lem.max.roots}},
\endalignedat
\]
resulting in $\bnd(\orb)\le\min\{10\cdot6,60\cdot1\}=60$.
\item\label{i.D6.3}
$\orb|_1\subset\cl2$, $\orb|_1\subset\cl3$, $\orb_3=\{0\}$,
$\orb|_4\subset\cl1$:
we have, for $k=2,4$,
\[*
\alignedat3
\cnt_1(\orb)&=1,&\quad\bnd_1(\orb)&=1&
 \quad&\text{by \autoref{D.2d} and}\\
\cnt_k(\orb)&=32,&\quad\bnd_k(\orb)&\le\bbD(6)=3&
 \quad&\text{by \autoref{D.1,3}\iref{D.1,3.1u}},
\endalignedat
\]
resulting in $\bnd(\orb)\le1\cdot3\cdot32=96$.
\item\label{i.D6.4}
$\orbk\subset\cl2$ for $k=1,\ldots,4$.
by \autoref{D.2d}, we have $\cnt_1(\orb)=\bnd_1(\orb)=1$ and, for $k=2,3,4$,
\[*
\cnt_k(\orb)=12,\quad\bnd_k(\orb)\le2,
\]
resulting in $\bnd(\orb)\le1\cdot2\cdot12\cdot12=288$.
However, since the last three blocks are special,
%\autoref{lem.D.special} and~
\eqref{eq.D.special} reduces this down to $1+1+12+12^2=158$.
\endroster
Summarizing, we arrive at
$\bnd(\Orb_h)\le8+3\cdot60+3\cdot96+288=764$, which is not quite
satisfactory. However, \autoref{lem.D.special} and~\eqref{eq.D.special}
reduce this bound down to $634$.

\section{Lattices with few roots: brute force}\label{S.few}

In this section, we consider the remaining Niemeier lattices rationally
generated by roots. This time, even the exact bounds $\bnd(\orb)$ for all
orbits result in $\bnd(\Orb_h)\ge720$,
and we need to deal with geometric rather than just admissible sets.

\theorem[see \autoref{proof.few}]\label{th.few}
Let $\bN:=\N(D)$ be a Niemeier lattice with $D$ one of
\[*
6\bD_4,\ 6\bA_4,\ 8\bA_3,\ 12\bA_2,\ \text{\rm or}\ 24\bA_1.
\]
Then, with five exceptions, \viz.~\eqref{eq.8A3.728} and
\eqref{eq.24A1.800.1}--\eqref{eq.24A1.728} below, one has $\ls|\fC|<720$ for
any $4$-polarization $h\in\bN$ and any geometric subset $\fC\subset\Orb_h$.
\endtheorem

\subsection{Proof of \autoref{th.few}}\label{proof.few}
We start with computing the sharp bounds $\bnd(\orb)$ for all combinatorial
orbits~$\orb$. If the additive bound $\bnd(\Orb_h):=\sum_\orb\bnd(\orb)$
exceeds $718$, we proceed by improving the similar bounds for (some) orbits
$\borb$ or unions~$\sS$ thereof.

Roughly, we try to prove that $\Bnd_m(\sS)=\varnothing$, see \eqref{eq.Bnd},
implying a reduced bound $\bb(\sS)\le m-2$.
Usually, the equality $\Bnd_m(\sS)=\varnothing$ above holds only modulo a number
of sets of large rank, which are singled out and analyzed in the course
of the proof (see mostly \tabrefs{\autoref*{s.rank} and \autoref*{s.push}});
it is these sets that give rise to the exceptions in \autoref{th.few} and,
eventually, the quartics in \autoref{th.main}.

The computation is done by the \emph{puzzle assembly} \tabref{s.puzzle}:
instead of adding individual conics one by one, we try to put together as
many and as large ``pieces'' as possible to fit a large number of conics into
a limited rank. The r\^{o}le of pieces of the puzzle is played by the elements
of the pre-computed set $\Bnd_0(\orb)$, $\orb\subset\borb$.

To ensure the convergence, we use a number of technical tricks:
\roster*
\item
\emph{clusters} $\orb\subset\cluster\subset\borb$ similar to homogeneous
chains in \autoref{s.chains}, see \tabref{s.clusters};
\item
computation \emph{up to rank $18$}, or even $17$ or~$16$, see
\tabref{s.rank}: if $\Bnd_m(\sS)$ consists of but a few sets of large rank,
these sets are analyzed by brute force;
\item
\emph{rank pushing}, see \tabref{s.push},
\ie, switching back to brute force for subsets
$\fC\subset\sS$ of rank close to $\rank\sS$.
\endroster
Further details are found in \cite{degt:800.tab}, where we also present the
topmost layer of the code used for each pair $h\in\bN$.
In the rest of this
section, we only show the five exceptional sets and a few otherwise interesting
examples: they appear in $\N(8\bA_3)$, see \autoref{s.ex.8A3},
and $\N(24\bA_1)$, see \autoref{s.ex.24A1}.
\qed

\convention\label{conv.ort}
In this and next section,
a saturated subset $\fC\subset\Orb_h$ is described by means
of its saturation $S:=\sat\fC$ which, in turn, is determined by the
polarized sublattice
\[*
h\in\ort\fC:=\bigl(h^\perp_S)^\perp_\bN\subset\bN.
\]
The reason for considering $\ort\fC$ instead of $S^\perp_\bN$ is purely
aesthetical: typically, it is generated by shorter vectors.
This is particularly useful in \autoref{S.Leech} below, where, dealing with the
Leech lattice~$\Leech$, it suffices to consider square~$4$ vectors only.
\endconvention

The original set $\fC$ is recovered as the set of vectors $l\in\bN$ such that
\[
l^2=4,\quad
2l\cdot v=h\cdot v\ \ \text{for each $v\in\ort\fC$}.
\label{eq.ort}
\]

\subsection{The lattice $\N(8\bA\sb3)$}\label{s.ex.8A3}
The lattice~$\bN$ is the index $256$ extension of $8\bA_1$ by the
discriminant classes
\[*
\cl3\oplus(\cl2\oplus\cl0\oplus\cl0\oplus\cl1\oplus\cl0\oplus\cl1\oplus\cl1),
\]
where the parenthesized expression runs over all seven cyclic permutations of its
arguments (see~\cite{Conway.Sloane}).
In particular, the kernel contains the class $\bigoplus_8\cl2$.

There are four $\OG(\bN)$-orbits of square~$4$ vectors $h\in\bN$; only one
of them, \viz. the one represented by
$h=h|_1\in\orbit4$ in $D_1$, results in large geometric subsets.

For the only exceptional set $\fC$ of cardinality~$728$,
the lattice $\ort\fC$ is generated by any three roots generating~$D_1$ and
two vectors $u=\bigoplus_ku|_k$, $v=\bigoplus_kv|_k$, where
\[
\label{eq.8A3.728}
%\gathered
%\ort\fC=\operatorname{Span}\{h,e_1,e_2,v_1,v_2\},\quad
% {\textstyle v_i=\bigoplus_kv_{i,k}},\quad
% v_{i,k}\in\cl2,\\
u|_k,v|_k\in\cl2,\quad
u|_{1}=v|_{1},\quad
u|_{1}\cdot h|_1=0,\quad
u|_{k}\cdot v|_{k}=0\quad\text{for $k=2,\ldots,8$}.
%\endgathered
\]
In the notation of \autoref{s.HAD}, assuming that $h|_1=\vv{\IS}$ and
letting $r:=\{1,2\}$, $s:=\{1,3\}$, we can take $u|_1=v|_1=v|_k=\be{r}$ and
$u|_k=\be{s}$ for $k=2,\ldots,8$.

%denoting by $e_1,e_2,e_3$ a standard root basis for~$\bA_3$ and
%assuming that $h=h|_1=e_1+2e_2+e_3$, we can take $u_1=v_k=\frac12(e_1+e_3)$
%for all $k=1,\ldots,8$ and $u_k=\frac12(e_1+2e_2+e_3)$ for $k=2,\ldots,8$.

\subsection{The lattice $\N(24\bA\sb1)$}\label{s.ex.24A1}
The lattice~$\bN$ is the index $4096$ extension of $24\bA_1$ generated by the
discriminant classes $\bigoplus_{k\in o}\cl1_k$, where $o\subset\Omega$ is a codeword in
the extended binary Golay code~$\Golay$, see,
\eg,~\cite{Conway.Sloane}.
There are two $\OG(\bN)$-orbits of square~$4$ vectors~$h$: either $h=\vv{s}$,
$s\subset\Omega$, $\ls|s|=2$, or $h=\frac12\vv{o}$, where $o\in\Golay$ is an
octad. In the latter case, there are no exceptional or otherwise interesting
geometric sets; we refer to \tabref{s.24A1-2} for the complete list of the
examples found.

Thus, denote $s:=\{1,2\}$ and assume that $h=\vv{s}$. Then, the conics are either
$l:=e_i\pm e_j$, $i\in s$, $j\in\Omega\sminus s$, or
\[
l:=\frac12\sum_{k\in o}\epsilon_ke_k,\quad
 \text{where $o\supset s$ is an octad, $\epsilon_1=\epsilon_2=1$,
 $\epsilon_k=\pm1$ for $k>2$}.
\label{eq.24A1.conics}
\]
Among the examples found, nine are of special interest; see \tabref{s.24A1-1}
for the complete list in terms of the original codewords as in \texttt{sage}.

\convention\label{conv.24}
We represent the generators $v_n:=\sum_kv_{n,k}e_k$ (the first always
being~$h$) of the lattice $\ort\fC$, see \autoref{conv.ort},
by pictographs, using the following notation for the
coefficients~$v_{n,k}$:
\minilist
\[*
({\black})\mapsto1,\quad
({\white})\mapsto-1,\quad
({\minussign})\mapsto\tfrac12,\quad
({\eqsign})\mapsto-\tfrac12,\quad
({\plussign})\mapsto\tfrac32.
\]
\endminilist
Needless to say that, for each~$n$, the set
$\{k\in\Omega\,|\,v_{n,k}\ne0\bmod\Z\}$ is a codeword.
For better transparency, we permute the index set~$\set$ so that
the subsets $\{1,\ldots,8\}$
and $\{7\ldots,14\}$ are among the octads.
Only the first $16$ positions are shown; it is understood that the last
character extends to the rest of the index set.
\endconvention

Given the first two generators,
the configuration is easily recovered as the set of
conics~\eqref{eq.24A1.conics} satisfying condition~\eqref{eq.ort}
\emph{for the last three generators~$v$} of $\ort\fC$. For each configuration~$\fC$
%of conics,
we show also the transcendental lattice(s) $T:=\NS(X)^\perp$ of the
quartic(s)~$X$ where this configuration is realized, in the inline notation
\[
[a,b,c]\quad\text{stands for}\quad
\Z u+\Z v,\quad u^2=a,\ u\cdot v=b,\ v^2=c,
\label{eq.inline}
\]
%$[a,b,c]$ for the rank~$2$ lattice $\Z u+\Z v$ with $u^2=a$, $u\cdot v=b$,
%$v^2=c$,
and the numbers $r+2c$ of the connected components of the
equiconical stratum (itemized according to~$T$), in the form $(r,c)$, where
\roster*
\item
$r$ is the number of real components and
\item
$c$ is the number of pairs of complex conjugate ones.
\endroster
These data are computed using Nikulin~\cite{Nikulin:forms} and,
for~\eqref{eq.24A1.608}, Miranda--Morrison~\cite{Miranda.Morrison:book},
similar to \autoref{proof.main} below.
%so that the total number of components is $r+2c$.
%, in the form $(r,c)$, the numbers~$r$ of real connected
%components and $c$ of pairs of complex conjugate components of the
\minilist
\chardef\minimax=17\relax
\[\label{eq.24A1.800.1}
\setsize{800}
\miniset
**......................\\
*.......................\\
..*.....................\\
--------................\\
------------------------\\
\endminiset
\T{
[4,0,40]&(1,0)\cr
}
\sym{[ 960, 11357 ]}
\]
\[\label{eq.24A1.800.2}
\setsize{800}
\miniset
**......................\\
*.......................\\
..*.....................\\
...*....................\\
------------------------\\
\endminiset
\T{
[4,0,40]&(1,0)\cr
}
\sym{[ 960, 11357 ]}
\]
\[\label{eq.24A1.736}
\setsize{736}
\miniset
**......................\\
*.......................\\
...*....................\\
--=-----................\\
------------------------\\
\endminiset
\T{
[12,0,16]&(1,0)\cr
}
\sym{[ 192, 1023 ]}
\]
\[\label{eq.24A1.728}
\setsize{728}
\miniset
**......................\\
*.......................\\
......*.................\\
......--------..........\\
------------------------\\
\endminiset
\T{
[14,0,14]&(1,0)\cr
}
\sym{[ 168, 42 ]}
\]
The first four sets are among the exceptions mentioned in \autoref{th.few}.
\[\label{eq.24A1.680}
\setsize{680}
\miniset
**......................\\
*.......................\\
........*...............\\
--=-----................\\
------------------------\\
\endminiset
\T{
[4,2,56]&(1,0)\cr
[16,6,16]&(1,0)\cr
}
\sym{[ 60, 5 ]}
\]
This graph is also realized with reducible conics, namely, by $X_{60}''$ in
\cite{DIS}: it has $60$ lines and $170+510=680$ conics, see
\autoref{rem.reducible} below.
\[\label{eq.24A1.660}
\setsize{660}
\miniset
**......................\\
*.......................\\
..*.....................\\
--------................\\
--------+---------------\\
\endminiset
\T{
[4,0,60]&(1,0)\cr
}
\sym{[ 60, 5 ]}
\]
This is the same abstract $K3$-surface as the Barth--Bauer quartic with $664$ conics,
see \cite{degt:4Kummer};
the latter is also embeddable to $\2\N(24\bA_1)$, see \tabref{s.24A1-1}.
\[\label{eq.24A1.656}
\setsize{656}
\miniset
**......................\\
*.......................\\
......*.................\\
......-=------..........\\
------------------------\\
\endminiset
\T{
[2,0,122]&(1,0)\cr
[10,4,26]&(0,1)\cr
}
\sym{[ 24, 12 ]}
\]
This is the largest known configuration of real conics on
a real quartic, see \autoref{proof.real}.
\[\label{eq.24A1.640}
\setsize{640}
\miniset
**......................\\
*.......................\\
..*.....................\\
---=----................\\
---+--------------------\\
\endminiset
\T{
[4,0,72]&(1,0)\cr
}
\sym{[ 192, 1023 ]}
\]
This is a Barth--Bauer $4^2\frak{A}_4$-quartics, \cf. \autoref{ad.family},
\eqref{eq.24A1.608}, and \autoref{proof.family}.
\[\label{eq.24A1.608}
\setsize{608}
\miniset
**......................\\
*.......................\\
..*.....................\\
...*....................\\
--------................\\
------------------------\\
\endminiset
\T{
{\rank}=19&(1,0)\cr
}
\sym{[ 192, 1023 ]}
\]
This is the largest known configuration of Picard rank~$19$, see
\autoref{proof.family}.
\endminilist

\section{The Leech lattice}\label{S.Leech}

The only Niemeier lattice left to be considered is the root free Leech
lattice~$\Leech$. Recall that there is a single $\OG(\Leech)$-orbit of
square~$4$ vectors $h\in\Leech$, see, \eg.,~\cite{Conway.Sloane}.

\theorem[see \autoref{proof.Leech}]\label{th.Leech}
Let $\Leech\ni h$ be a $4$-polarized Leech lattice. Then,
with four exceptions~\eqref{eq.Leech.800}--\eqref{eq.Leech.720} below, one has
$\ls|\fC|<720$ for any geometric set $\fC\subset\Orb_h$.
\endtheorem

\subsection{Iterated index~$2$ subgroups}\label{s.subgroups}
There is a single orbit $\borb=\Orb_h$, whereas, since $\Leech$ is root free,
we have $\RG(\Leech)=1$ and all combinatorial orbits are singletons. These
facts make the application of the puzzle assembly algorithm problematic.
Besides, since there are no roots, almost any subset $\fC\subset\Orb_h$ is
admissible, see \autoref{def.admissible},
whereas it is the inadmissibility that
rules out most intermediate sets in the other Niemeier lattices.

Therefore, we have to shift the paradigm and, instead of starting
from~$\varnothing$ and working upwards, we start from~$\Orb_h$ and go
downwards, searching for large \emph{geometric} subsets $\fC\subset\Orb_h$. In
fact, most of the time the principal criterion is $\rank\fC\le20$.

The following lemma has essentially appeared in~\cite{DIO}. We present its
formal proof since it is to be followed literally by the new version of the
algorithm.

\lemma\label{lem.subgroups}
Let $\bN\ni h$ be a $4$-polarized Niemeier lattice. Then, for any saturated subset
$\fC\subset\Orb_h$, there is a chain
\[*
\fC=\fC_n\subset\fC_{n-1}\subset\ldots\subset\fC_1\subset\fC_0=\Orb_h
\]
such that, for each $k=1,\ldots,n$, one has $\fC_{k}=\Orb_h\cap S_{k}$
for a certain polarized index~$2$ sublattice $S_{k}\subset\spnZ\fC_{k-1}$,
$S_k\ni h$.
\endlemma

\proof
Clearly, since $\spnZ\Orb_h$ is generated by conics, there is a sequence
\[*
\fC=\fS_m\subset\fS_{m-1}\subset\ldots\subset\fS_1\subset\fS_0=\Orb_h
\]
of \emph{saturated} subsets such that $\rank\fS_k=\rank\fS_{k-1}-1$ for each
$k=1,\ldots,m$. Fix a pair $\fS_{k+1}\subset\fS_{k}$ of consecutive sets,
let $\fS_{k,0}:=\fS_k$, and denote by
\[*
p_0\:\spnZ\fS_{k,0}\to Q_0:=\spnZ\fS_{k,0}/\spnZ\fS_{k+1}
\]
the quotient projection. The group $Q_0\cong\Z$ is generated by (the images
of) finitely many conics $l\in\fS_{k,0}\sminus\fS_{k+1}$, and it is immediate
that the cardinality of the set
\[*
\fS_{k,1}:=\fS_{k,0}\cap p_0\1(2Q_0)\supset\fS_{k+1}
\]
is strictly less than $\ls|\fS_{k,0}|$. Iterating this procedure, we arrive
at a chain
\[*
\fS_{k}=\fS_{k,0}\supset\fS_{k,1}\supset\ldots\supset\fS_{k,r}=\fS_{k+1}
\]
satisfying the ``index $2$ sublattice'' condition; it is bound to terminate
as all inclusions are proper and $\ls|\fS_k\sminus\fS_{k+1}|<\infty$.
%and this sequence must terminate: $\fS_{k-1,r}=\fS_k$.
Replacing each pair $\fS_{k+1}\subset\fS_k$ of consecutive
sets with such a chain, we obtain a sequence as in the statement.
\endproof

Note that we do not assert that whenever the index $[\spnZ\fC_p:\spnZ\fC_q]$
is finite it is a power of~$2$. Though, in our computation it was always the
case: before the rank drops, we observed indices up to~$32$.

\autoref{lem.subgroups} is used to construct large geometric subsets
of~$\Orb_h$. We start with $\Orb_h$ (or, in some cases, with another
large admissible saturated set) and construct chains
\[*
\Orb_h=:\fC_0\supset\fC_1\supset\ldots
\]
recursively: once a set $\fC_k$ has been constructed, we
\roster
\item
consider the projection
$(\bmod\,2)\:\spnZ\fC_k\to V_k:=(\spnZ\fC_k)\otimes\FF$;
\goodbreak
\item
compute the $(\stab\fC_k)$-orbits on the annihilator
$h^\perp\subset V_k\dual$;
\goodbreak
\item
for a representative $v$ of each orbit, let
$\fC_{k+1}:=\fC_k\cap(\bmod\,2)\1(v^\perp)$;
\goodbreak
\item
discard the subsets~$\fC_{k+1}$ of cardinality less than a preset threshold.
\endroster
Unfortunately, this procedure, used in~\cite{degt:conics,DIO}, diverges at a
rate unacceptable for the current problem. To make it more tame and reduce
the overcounting, we follow the proof of \autoref{lem.subgroups} more
literally. More precisely, at each step~$\fS_k$, we
\roster*
\item
keep track of the ``original'' saturated set $\bar\fC_k=\sat\fC_k$;
\item
discard~$\fC_{k+1}$ if $\rank\fC_{k+1}=\rank\fC_k$ and
$\spnZ\bar\fC_k/\spnZ\fC_{k+1}$ is not cyclic;
\item
discard~$\fC_{k+1}$ if $\rank\fC_{k+1}<\rank\fC_k$ and $\fC_{k+1}$ is not
saturated.
\endroster
Besides, we perform the computation rank by rank: once all saturated sets of
some rank~$r$ have been collected, before proceeding we leave a single
representative of each $\OG_h(\bN)$-orbit only. In fact, since we are
interested in the existence of large subsets only (not in a particular
embedding to~$\Orb_h$) and the computation depends on $\spnZ\bar\fC_k$,
we
go one step further and
leave a single representative of each graph isomorphism class.
Other technical details are discussed in~\tabref{S.Leech.tab}.

\subsection{Proof ot \autoref{th.Leech}}\label{proof.Leech}
The lattice $\Leech\subset\bH_{24}\bigl(\frac18\bigr)$ is generated by the
square $4$ vectors of the form
\roster*
\item
$4\cdot\vv{u}$, where $u\subset\IS$ and $\ls|u|=2$,
\item
$2\cdot\vv{o}$, where $o\subset\IS$ is an octad of the Golay code, or
\item
$\vv{\IS}-4e_i$, where $i\in\IS$,
\endroster
as well as all vectors obtained from these by the simultaneous reversal of
the sign at the elements of a codeword of the
extended binary
Golay code. This
description does not reflect the full automorphism
group $\OG(\Leech)$: the latter is transitive on the
square~$4$ vectors $h\in\Leech$ (see, \eg, \cite{Conway.Sloane}).

Thus, we fix a $4$-polarization $h\in\Leech$ and
start with the $h$-even index~$2$ sublattice
$\bar\Leech:=\spnZ\Orb_h\subset\Leech$;
it has larger orthogonal group $\bigl[\OG_h(\bar\Leech):\OG_h(\Leech)\bigr]=2$.
Besides, since $\Leech$ is root free, the set $\Orb_h$ itself and any subset
thereof is admissible in $\bar\Leech$.
It remains to invoke \autoref{lem.subgroups} and list all (graph isomorphism
classes of) subsets $\fC\subset\Orb_h$ of rank $\rank\fC\le20$,
see~\eqref{eq.rank}, and size $\ls|\fC|\ge720$. With the new improvements
(see \autoref{s.subgroups} and
\tabrefs{\autoref*{s.Leech.literal}--\autoref*{s.Leech.failures}})
this takes less than two
days (\vs. two months in the version of~\cite{DIO}) and we arrive at nine
sets~$\fC$:
one has $\rank\fC=20$ and
\[*
\ls|\fC|=896,
\hyperref[eq.Leech.800]{\ul{800}},768,760,740,
\hyperref[eq.Leech.736]{\ul{736}},736,
\hyperref[eq.Leech.728]{\ul{728}},
\hyperref[eq.Leech.720]{\ul{720}}.
\]

Strictly speaking, found are graph isomorphism classes rather than
$\OG_h(\Leech)$-orbits; hence, we continue the analysis in terms of the
lattice $\spnZ\fC$, which is a graph invariant.
For each set~$\fC$, we use
Nikulin's theory~\cite{Nikulin:forms} to compute the genus
%$\operatorname{genus}(\ort\fC)$
of the lattice $\ort\fC$ (see \autoref{conv.ort})
and, referring to Gordon L. Nipp's
tables~\cite{Nipp:quaternary.forms,Nebe.Sloane:lattices}, we find that there
is but a single
\emph{root free} representative $\ort\fC$. Furthermore, all classes in the
genus of any proper finite index extension do have roots and, hence, cannot
be embedded to~$\Leech$.

We conclude that the sublattice $\spnZ\fC\subset\Leech$ must be primitive
and, hence, it can be used in \autoref{def.geometric}.
Only the four underlined sets are geometric; the corresponding lattices
$\ort\fC$ are shown below, both as the Gram matrix and as a quintuple of
square~$4$ generators.
(The remaining five sets are
found in \tabrefs{(\ref*{eq.Leech.nongeometric})}.)
For the generators, we adopt \autoref{conv.24}, except
that we use the notation
\[*
({\black})\mapsto4,\quad
({\white})\mapsto-4,\quad
({\plussign})\mapsto2,\quad
({\minussign})\mapsto-2
\]
for the coefficients, assume that $\{1,\ldots,8\}$ and $\{5,\ldots,12\}$ are
octads, and cut the display at position~$14$, extending each vector by
zeroes.
\bgroup
\let\ministrut\strut
\[
\ls|\fC|=800\:\
\bmatrix
4&0&0&0&2\\
0&4&0&1&2\\
0&0&4&1&2\\
0&1&1&4&0\\
2&2&2&0&4\endbmatrix\quad
\miniset
....*..!.....\\
.....*..*....\\
....++-+-+++.\\
++++++++.....\\
....++--++++.\\
\endminiset
\label{eq.Leech.800}
\]
\[
\ls|\fC|=736\:\
\bmatrix
4&2&0&0&0\\
2&4&0&0&2\\
0&0&4&1&2\\
0&0&1&4&2\\
0&2&2&2&4\endbmatrix\quad
\miniset
!...*........\\
....**.......\\
......*.*....\\
++++++++.....\\
....++++++++.\\
\endminiset
\label{eq.Leech.736}
\]
\[
\ls|\fC|=728\:\
\bmatrix
4&2&2&2&2\\
2&4&0&1&0\\
2&0&4&0&1\\
2&1&0&4&2\\
2&0&1&2&4\endbmatrix\quad
\miniset
....**.......\\
++++++++.....\\
+---++-+.....\\
....+++--+++.\\
....++--++++.\\
\endminiset
\label{eq.Leech.728}
\]
\[
\ls|\fC|=720\:\
\bmatrix
4&0&0&2&2\\
0&4&0&0&2\\
0&0&4&0&2\\
2&0&0&4&1\\
2&2&2&1&4\endbmatrix\quad
\miniset
**...........\\
..**.........\\
....**.......\\
*.......*....\\
++++++++.....\\
\endminiset
\label{eq.Leech.720}
\]
%\end{gather}
\egroup
Since only square~$4$ vectors are involved, it is not difficult to show, by
brute force, that each of the four lattices $\ort\fC$ above admits a unique, up to
$\OG_h(\Leech)$, polarized isometry $(\ort\fC\ni h)\into(\Leech\ni h)$,
completing the proof of the uniqueness.
\qed

\section{Proofs of the main results}\label{S.proofs}
In this section we collect our previous findings and fill in the missing
parts of the proofs of the principal results stated in the introduction.

\subsection{Proof of \autoref{th.main}}\label{proof.main}
According to \autoref{prop.reduction}, the graph of conics of a smooth
quartic is isomorphic to a saturated admissible set~$\fC$ in a $4$-polarized
Niemeier lattice and, by Theorems~\ref{th.many}, \ref{th.few},
and~\ref{th.Leech}, up to the relevant orthogonal groups,
there are but nine such sets of cardinality $\ls|\fC|\ge720$, \viz.
\[*
\eqref{eq.24A1.800.1}\cong\eqref{eq.24A1.800.2}\cong\eqref{eq.Leech.800},\quad
\eqref{eq.24A1.736}\cong\eqref{eq.Leech.736},\quad
\eqref{eq.8A3.728}\cong\eqref{eq.24A1.728}\cong\eqref{eq.Leech.728},\quad
\eqref{eq.Leech.720};
\]
here, the abstract graph isomorphisms are established using the
\texttt{digraph} package in \GAP~\cite{GAP4.13}.
For each of the nine sets~$\fC$ we have $S:=\spn\fC=\spnZ\fC$ and, hence, the
rest of the computation depends on the abstract graph isomorphism class only,
leaving but four cases.

For each of the first three graphs,
%using Nikulin's theory~\cite{Nikulin:forms}, we conclude that
the modified lattice $(S\ni h)\smod$, see \autoref{s.S*},
admits a unique (up
to the group $\OG^+(\bL)$ of auto-isometries of~$\bL$
preserving the orientation of maximal
positive definite subspaces) primitive isometry $S\smod\!\into\bL$, resulting
in a single projective equivalence class of quartics (see, \eg,
Dolgachev~\cite{Dolgachev:polarized}). On the other hand,
%by~\cite{Nikulin:forms} again,
$(S\ni h)\smod$ has no
$h$-odd index~$2$ extensions; hence, the configuration is realized by
line-free quartics only.
For both assertions, the computation employing
Nikulin's discriminant forms~\cite{Nikulin:forms}
is considered quite straightforward nowadays and
therefore is left to the reader. We have it implemented in \GAP.

For~\eqref{eq.Leech.720}, the situation is the opposite: there is no
primitive isometry $S\smod\into\bL$ but, up to $\OG_h(S\smod)$, there is a
unique $h$-odd index~$2$ extension $S'\supset S\smod$, and the latter does
admit a unique primitive isometry $S'\into\bL$, resulting in a single
quartic~$X$. Upon computing lines and conics on~$X$ by means of~\eqref{eq.conics},
we identify it as Schur's quartic (as the only smooth quartic
with $64$ lines, see~\cite{DIS}). Alternatively, Schur's is the only quartic
with the transcendental lattice $[8,4,8]$,
see~\eqref{eq.inline} and~\cite{degt:singular.K3}.
\qed

\remark\label{rem.Aut}
For each quartic~$X$ encountered in this paper, including the four in
\autoref{th.main}, the lattice $\NS(X)$ is generated over~$\Z$ by lines and
conics; hence, the group $\OG_h(\NS(X))=\Aut(\Fn_*X)$ can easily be computed
using the \texttt{digraph} package in \GAP~\cite{GAP4.13}. Then one can also
compute the groups $\Sym_hX\subset\Aut_hX$ of (symplectic)
projective automorphisms of~$X$. In particular, $\Sym_hX=M_{20}$ for
quartic~\iref{i.M20} in \autoref{th.main} (first observed by X.~Roulleau,
private communication) and, due to the
uniqueness~\cite{Bonnafe.Sarti,degt:4Kummer} of a quartic with a symplectic
action of~$M_{20}$, one can reuse the equation found in~\cite{Mukai}.
\endremark

\subsection{Proof of \autoref{ad.real}}\label{proof.real}
As explained in~\cite{DIS} (in the context of lines, but the argument extends
to rational curves of any degree), when estimating the maximal number of real
conics on a real quartic it suffices to assume that \emph{all} conics on
a quartic~$X$
are real and that they generate $\NS(X)$ over~$\Q$. Furthermore, if $\NS(X)$
is rationally generated by conics, then $X$ admits a real structure with all
conics real if and only if the transcendental lattice $T:=\NS(X)^\perp$
contains $\Z a$, $a^2=2$, or $\bU(2)$ as a sublattice.
Thus, \autoref{ad.real} is an immediate consequence of an analysis of the
examples found in the course of the proof; the quartic is
\eqref{eq.24A1.656}.
%in $\N(24\bA_1)$.
\qed

\subsection{Proof of \autoref{ad.family}}\label{proof.family}
The configuration of conics is \eqref{eq.24A1.608} in $\N(24\bA_1)$,
which is of rank~$19$. The transcendental lattice $T:=\NS(X)^\perp$ of a
generic member $X\in\Cal{X}$ is
\[*
T\cong[-4]\oplus[8]\oplus[8]\cong[4]\oplus[-8]\oplus[-8],
\]
where $[s]$ stands for the rank~$1$ lattice $\Z a$, $a^2=s$. The group
$\Sym_hX\cong4^2\frak{A}_4$ is computed in terms of
the Fano graph $\fC=\Fn X$, see
\autoref{rem.Aut}. The geometric oversets $\fC'\supset\fC$ are found by
applying \autoref{lem.subgroups} (with all the improvements to the algorithm
discussed in
\tabrefs{\autoref*{s.Leech.literal}--\autoref*{s.Leech.failures}}) to all
other examples found in the course of the proof
for this particular polarization of $\N(24\bA_1)$; they are
\eqref{eq.24A1.800.1}, \eqref{eq.24A1.800.2},
\eqref{eq.24A1.736}, and~\eqref{eq.24A1.640}.
Then, general theory of $K3$-surfaces implies that any geometric overset can be
realized by a degeneration of quartics.
\qed

\remark\label{rem.rank19}
Applying \autoref{lem.subgroups} to the other known examples, we could not
find a rank~$19$ geometric set larger than~\eqref{eq.24A1.608} in
\autoref{ad.family}.
\endremark

\remark\label{rem.oversets}
We do not assert that the list of (graph isomorphism classes of) oversets in
the proof of \autoref{ad.family} is complete. Probably, the complete list of
degenerations can be found by analyzing corank~$1$ abstract graph extensions
of~$\fC$ in~\eqref{eq.24A1.608}, but we leave this exercise to the reader.
\endremark

\subsection{Concluding remarks}\label{s.remarks}
A great deal of open questions related to conics on quartic surfaces are left
beyond the scope of this paper. Thus, we do not discuss fields of definition
of positive characteristic or quartics with $\bA$--$\bD$--$\bE$
singularities, which are still $K3$-surfaces. (We did try to apply an
analogue of \autoref{lem.subgroups} to the set $\Orb_h\subset\N(24\bA_1)$,
but we could not find a configuration with more than a couple of hundreds of
conics. Note that, in the presence of singularities, the total number of \emph{all} conics,
including those containing exceptional divisors as components, may be much
larger than $800$.)

It is not immediately clear whether the examples given by
Addenda~\ref{ad.real} and~\ref{ad.family} are, indeed, maximal.
Besides, even though \autoref{th.main} does give us the bound $\Bs(4,2)=720$,
it is not clear how the maximal number of \emph{irreducible} conics may be
affected by the presence of lines or exceptional divisors. In all examples in
the next remark, in the presence of \emph{many} lines, most conics are
reducible.

\remark\label{rem.reducible}
According to~\cite{degt:4Kummer}, the maximal number of \emph{reducible}
conics on a smooth quartic is $576$, attained at Schur's classical
quartic~\iref{i.T192} in \autoref{th.main}.
As an example, we computed the configurations of conics on all smooth
quartics with more than $48$ lines (see \cite{DIS,degt.Rams:quartics}).
Only six of them have more than $600$
conics, \viz. precisely those with at least $56$ lines:
\begin{alignat*}3
X_{64}\:&\quad\text{$64$ lines},&\quad&\text{$144 + 576 = 720$ conics},
 &\quad&\text{no planes, see \eqref{eq.Leech.720}},\\
X_{60}'\:&\quad\text{$60$ lines},&\quad&\text{$140 + 500 = 640$ conics},
 &\quad&\text{no planes, see \tabref{s.12A2-2}},\\
X_{60}''\:&\quad\text{$60$ lines},&\quad&\text{$170 + 510 = 680$ conics},
 &\quad&\text{$10$ planes, see \eqref{eq.24A1.680}},\\
X_{56}\:&\quad\text{$56$ lines},&\quad&\text{$184 + 440 = 624$ conics},
 &\quad&\text{$16$ planes, see \tabref{s.24A1-2}},\\
Y_{56}\:&\quad\text{$56$ lines},&\quad&\text{$188 + 448 = 636$ conics},
 &\quad&\text{$20$ planes},\\
Q_{56}\:&\quad\text{$56$ lines},&\quad&\text{$208 + 448 = 656$ conics},
 &\quad&\text{$24$ planes, see \tabref{s.24A1-2}}.
\end{alignat*}
Here, by a \emph{plane} we mean a plane section split into two irreducible
conics: formally, these conics are undetectable in terms of lines only.

Remarkably, most configurations appear on our list of examples,
substantiating the suggestion that it must be close to complete, even though
some gaps do not look quite convincing.
\endremark

{
\let\.\DOTaccent
\def\cprime{$'$}
\bibliographystyle{amsplain}
\bibliography{degt}
}

\end{document}